%% file: main.tex
\documentclass[reqno,11pt]{amsart}
\usepackage{color}
\usepackage{comment}
\usepackage{pdfpages}
\usepackage{graphicx}
\usepackage{dcolumn}
\usepackage{dsfont}
\usepackage{tabularx}
\usepackage{yfonts}
\usepackage{tikz-cd}
\usetikzlibrary{shapes.misc}
\usetikzlibrary{shapes.symbols}
\usetikzlibrary{snakes}
\usetikzlibrary{decorations}
\usetikzlibrary{decorations.markings}

\newenvironment{chosenEnum}[1][]{%
\enumerate[label=(\roman*), #1]%
}
{
\endenumerate
}

\usepackage[
top    = 3.5cm,
bottom = 3.00cm,
left   = 2.8cm,
right  = 2.8cm]{geometry}

\RequirePackage[numbers]{natbib}
\RequirePackage[colorlinks=true, pdfstartview=FitV, linkcolor=blue,
  citecolor=blue, urlcolor=blue]{hyperref}
\usepackage{paralist}

\usepackage{tikz} 
\usepackage{dcolumn}
\usetikzlibrary{arrows,positioning}
\tikzset{
    >=stealth',
    pil/.style={
           ->,
           thick,
           shorten <=2pt,
           shorten >=2pt,}
}

\usepackage{graphicx}
\graphicspath{{./figures/}}
\usepackage{amsfonts}
\usepackage{amsmath}
\usepackage{amsthm}
\usepackage{amssymb}
\usepackage{amsbsy}
\usepackage{epsfig}
\usepackage{natbib, mathrsfs}
\usepackage{verbatim}
\usepackage[latin1]{inputenc}
\usepackage{mhequ}
\usepackage{algorithm}
\usepackage{algorithmic}

\usepackage{import}
\usepackage{enumitem} 
\usepackage{caption,subcaption}

        \let\l=\lambda

 \let\x=\xi

\renewcommand{\leq}{\;\leqslant\;}                   
\renewcommand{\geq}{\;\geqslant\;}                   

\newcommand\weakconv{\stackrel{D}{\longrightarrow}}

\newcommand\weakeq{\stackrel{D}{=}}

\newcommand\R{\mathbb{R}}
\newcommand\E{\mathbb{E}}

\newcommand\genFin{{{\mathcal{L}}_{\theta}}}
\newcommand\genInf{{\Hat{\mathcal{L}}}}
\newcommand\genPD{{\mathcal{G}_\theta}}
\newcommand\mutFin{{A_\theta}}
\newcommand\mutInf{{\Hat{A}}}

\newcommand{\N}{\mathbb N}

\definecolor{WowColor}{rgb}{.75,0,.75}
\definecolor{SubtleColor}{rgb}{0.9,0,0}

\newcounter{margincounter}

\newcounter{latercounter}

\newtheorem{theorem}{Theorem}[section]
\newtheorem{lemma}[theorem]{Lemma}
\newtheorem{proposition}[theorem]{Proposition}
\newtheorem{corollary}[theorem]{Corollary}
\newtheorem{remark}[theorem]{Remark}

\newtheorem*{question*}{Question}

\newtheorem*{remark*}{Remark}
\newtheorem*{idefinition*}{Definition}
\newtheorem*{example*}{Example}

\begin{document}

\title{Size-biased diffusion limits and the inclusion process}

\author[P. Chleboun]{Paul Chleboun}
\author[S. Gabriel]{Simon Gabriel}
\author[S. Grosskinsky]{Stefan Grosskinsky}
\address{P. Chleboun, Department of Statistics, University of Warwick,  Coventry,  CV4 7AL,  United Kingdom } 
\email{paul.i.chleboun@warwick.ac.uk}
\address{S. Gabriel, Mathematics Institute, University of Warwick, Coventry CV4
7AL, United Kingdom \\
and  Institute for Analysis and Numerics,
University of M\"unster, 48149 M\"unster, Germany} 
\email{simon.gabriel@uni-muenster.de}
\address{S. Grosskinsky, Department of Mathematics, University of Augsburg, 86159 Augsburg, Germany} 
\email{stefan.grosskinsky@math.uni-augsburg.de}
\thanks{
S. Gabriel thanks Tommaso Rosati for helpful discussions and pointing out the
reference \cite{DeMarco11}.
S. Gabriel is supported by the Warwick Mathematics Institute Centre for Doctoral Training, and acknowledges funding from the University of Warwick and EPSRC through grant EP/R513374/1.
}


\begin{abstract}
We consider the inclusion process on the complete graph with vanishing diffusivity, which leads to condensation of particles in the thermodynamic limit. Describing particle configurations in terms of size-biased and appropriately scaled empirical measures of mass distribution, we establish convergence in law of the inclusion process to a measure-valued Markov process on the space of probability measures. In the case where the diffusivity vanishes like the inverse of the system size, the derived scaling limit is equivalent to the well known Poisson-Dirichlet diffusion, offering an alternative viewpoint on these well-established dynamics. Moreover, our novel size-biased approach provides a  robust description of the dynamics, which covers all scaling regimes of the system parameters and yields a natural extension of the Poisson-Dirichlet diffusion to infinite mutation rate. We also discuss in detail connections to known results on related Fleming-Viot processes.
\end{abstract}

\maketitle


\setcounter{tocdepth}{1}
\tableofcontents

\import{.}{bulk.tex}

\appendix

\import{.}{appendix.tex}

\bibliography{cites}
\bibliographystyle{alpha}

\end{document}

%% file: bulk.tex
\section{Introduction and main results} 

The inclusion process $(\eta^{(L,N)}(t))_{t\geq 0}$ on the complete graph of
$L$ sites, describing the interaction of $N$ particles, is characterised by the
infinitesimal generator 
\begin{align}\label{eq_def_ip_dynamics}
    \mathfrak{L}_{L,N}f(\eta)
    :=
    \sum_{\substack{x,y=1 \\ x\neq y}}^L 
    \eta_x(d+\eta_y) [f(\eta^{x,y})-f(\eta)]\,.
\end{align}
Here $\eta\in \Omega_{L,N}:= \{\eta \in \N_0^L \,:\, \sum_{x=1}^L \eta_x=N\}$ and $f:\Omega_{L,N}\mapsto \R$. Moreover, $\eta^{x,y}$ denotes the configuration $\eta$ where one particle moved from $x$ to $y$, i.e. $\eta^{x,y} =\eta+e^y-e^x$, with $(e^{x})_{z}:= \delta_{x,z}$, provided $\eta_x>0$.
The constant $d\geq 0$ is usually referred to as \emph{diffusivity}.

The dynamics of the inclusion process consists of two parts. Each particle performs a continuous-time random walk with rate $d$, and in addition particles attract each other at unit rate. 
The process was introduced in \cite{GKR07} as the dual of a model of energy
transport and is the natural counterpart of the exclusion process, since interactions are attractive rather than repulsive.
Moreover, it can be interpreted in the context of population genetics,
describing a population of Moran type \cite{Mo58} where $d$ corresponds to the
mutation rate and resampling occurs at rate $1$. 

The present article concerns the study of the inclusion process in the thermodynamic limit: $N,L\to\infty$ such that $N/L\to \rho\geq 0$, usually abbreviated by simply writing $N/L\to \rho\geq 0$.
We rescale the diffusivity $d=d(L)$ with the size of the system such that $d\to 0$ as $L\to \infty$. This leads to clustering of particles into chunks of diverging size \cite{GRV11,GRV13}. In the context of stochastic particle systems this phenomenon is known as condensation. 
The occurrence of condensation and the statistics of the condensed phase have been
studied for various particle systems (see e.g. \cite{EH05,GL12,CG14}),
references for the inclusion process and related models include
\cite{WaEv12,CCG14,CCG15,CGG21} and \cite{BDG17,KS21} in the context of
metastable dynamics of a single condensate site. These results apply in different scaling regimes with fixed volume and diverging mass which are not addressed in this paper, and we discuss the connections to our results in detail in Section \ref{sec_overview}. In \cite{JCG19} the condensed phase in the
stationary inclusion process was shown to exhibit clusters on scale $L$ with a
Poisson-Dirichlet size distribution in the case $dL\to\theta\in\R_+$, whereas
for $dL\to\infty$ the clusters are of scale $1/d$ with independent exponential
distribution. The aim of this paper is to characterise the dynamics of the
mass distribution in the condensed phase in both regimes in terms of diffusion limits. 


\subsection{Main results}

Our main results are twofold. 
On the one hand, we describe convergence of the inclusion process towards measure-valued diffusions, which are parameterised by $\theta:= \lim_{N/L \to \rho} dL\in [0,\infty]$.
In the case $\theta<\infty$, the determined scaling limit is equivalent to the well known infinitely-many-neutral-alleles diffusion model \cite{EK81}, also known as Poisson-Dirichlet diffusion.
On the other hand, the new description of the infinitely-many-neutral-alleles
diffusion model, we give in the following, allows us to  construct a natural extension of the process, when $\theta=\infty$.

\subsubsection{Scaling limits of the inclusion process}

We distinguish between the cases $\theta<\infty$ and $\theta=\infty$.
When $\theta<\infty$, we embed particle configurations into $\mathcal{M}_1([0,1])$ using the maps
\begin{align}\label{eq_def_dLfin_embedding}
\mu^{(\cdot)}_{L,N}: \Omega_{L,N} \to \mathcal{M}_1([0,1])\quad\mbox{of the form}\quad    \mu^{(\eta)}_{L,N}:= \sum_{x=1}^L \frac{\eta_x}{N}\delta_{\frac{\eta_x}{N}}\,.
\end{align}
Here $\mathcal{M}_1([0,1])$ denotes the set of probability measures on
$[0,1]$, equipped with the topology induced by weak convergence of measures.
Note that this corresponds to a size-biased empirical measure on the space of mass fractions $[0,1]$. Not every measure in $\mathcal{M}_1([0,1])$ can be approximated by particle configurations using $\mu^{(\cdot)}_{L,N}$, e.g. every point mass $\alpha \delta_z$ in \eqref{eq_def_dLfin_embedding} satisfies $\alpha\geq z$.
Instead, we restrict ourselves to the closed subspace of atomic measures
\begin{equation}\label{eq_def_E}
E:=\mu^{(\overline{\nabla})}\subset \mathcal{M}_1([0,1])\ ,    
\end{equation}
defined as the range of the function
\begin{align}\label{eq_embedding_kingman_to_measure}
    \mu^{(\cdot)} : \overline{\nabla} \to \mathcal{M}_1([0,1])\quad\mbox{with}\quad\mu^{(p)} = (1-\|p\|_1)\delta_0+ \sum_{i=1}^\infty p_i \delta_{p_i} \,,
\end{align}
where $\overline{\nabla}$ denotes the Kingman simplex
\begin{equation}\label{eq_def_kingsim}
\overline{\nabla} =\{p\in [0,1]^{\N}\,:\, p_1\geq p_2\geq \ldots\,, \
\sum_{i=1}^\infty p_i\leq 1\}\ ,
\end{equation}
equipped with the product topology induced by $[0,1]^{\N}$.
We may often drop the subscripts in \eqref{eq_def_dLfin_embedding} and simply
write $\mu^{(\cdot)}$, but the meaning will be clear from the context. 
The above mappings $\mu^{(\cdot)}$ and $\mu^{(\cdot)}_{L,N}$ do not preserve spatial information of particle configurations, but this also does not enter the dynamics on the complete graph. Note that the map \eqref{eq_embedding_kingman_to_measure} was already mentioned in \cite{EK81}, however, only to prove denseness of the domain of functions considered there.

Moreover, any \emph{(unordered)} particle configuration $\eta\in \Omega_{L,N}$ can be mapped onto an
element in $\overline{\nabla}$ using the \emph{ranked mass embedding} 
\begin{equation}\label{eq_ranked_mass}
\begin{aligned}
\tfrac{1}{N}\Hat{\eta} := \tfrac{1}{N} ( \Hat\eta_{1},  \Hat\eta_{2}, \ldots,
\Hat\eta_{L})\,,
\end{aligned}
\end{equation} 
with $\Hat{\eta}$ representing the entries of $\eta$ in decreasing order.
Another representation of $ \eta $ is by \emph{size-biased sampling}.
More precisely, we define 
\begin{equation}\label{eq_def_sb}
\begin{aligned}
\Tilde{\eta}_k:=\eta_{\sigma(k)}\,,
\end{aligned}
\end{equation}
for some random permutation $\sigma$ of $\{ 1,\ldots ,L\}$ generated iteratively as follows: First
\begin{chosenEnum}
    \item $\sigma(1)=x$ with probability $\frac{\eta_x}{N}$, $x\in \{1,\ldots, L\}$\,,
\end{chosenEnum}
and for any following index $k=2, \ldots, L$
\begin{chosenEnum}
    \item[(ii)] $\sigma(k)=x$ with probability $\frac{\eta_x}{N-\sum_{j=1}^{k-1}\eta_{\sigma(j)}}$, $x\in \{1,\ldots, L\}\setminus \{\sigma(1), \ldots, \sigma(k-1)\}$\,.
\end{chosenEnum}
We refer to \cite[Definition 2]{JCG19} for details.
The concept of size-biased sampling plays a central role in the present
article, since \eqref{eq_def_dLfin_embedding} can be thought of as the empirical distribution of the first size-biased marginal $ \Tilde\eta_{1}$.\\


In order to describe the limiting dynamics, we consider the domain of functions
\begin{align}\label{eq_def_DL}
    \mathcal{D}({\genFin}) = \text{sub-algebra of } C(E) \text{ generated by functions } \mu \mapsto \mu(h)\,, \ h\in C^3([0,1])\,.
\end{align}
\sloppy The pre-generator of the corresponding superprocess acting on a function $H(\mu) = \mu(h_1)\cdots \mu(h_n)\in \mathcal{D}({\genFin})$ then reads
\begin{align}\label{eq_def_L}
    \genFin H(\mu) 
    &:= 
    2\sum_{1\leq k < l \leq n} \big( \mu(Bh_k Bh_l) -\mu(Bh_k) \mu(Bh_l)  \big)
    \prod_{j\neq k,l} \mu(h_j)\\
    &\qquad +\sum_{1\leq k \leq n} \mu({\mutFin} h_k)\prod_{j\neq k} \mu(h_j)\,, \nonumber
\end{align}
with $Bh(z) := h(z) +zh'(z)=(z h(z))'$.
Here, the first part is usually referred to as interaction term and ${\mutFin}$ denotes the single-particle operator of the form
\begin{align}\label{eq_def_A_singlePart}
    {\mutFin} h(z):=& 
     (1-z)(Bh)'(z)  +\theta(Bh(0)-Bh(z))\\
    =&
    z(1-z) h''(z) +(2-z(2+\theta ) )h'(z) +\theta(h(0)-h(z))
    \,.\nonumber
\end{align}
The operator $Bh(z)$ should be thought of as a `size-biased derivative', which appears due to our choice of embedding \eqref{eq_def_dLfin_embedding}.
For example, a single site containing a mass fraction $z\in [0,1]$ will be represented by a point-mass $z \delta_z$. Thus, change in $z$ will result both in a change of the amount of mass and its position. 

Our first result identifies the process described by $\genFin$ as the correct scaling limit of the inclusion process. Here and in the following we will sometimes write $\mu_{\#} \eta$ instead of $\mu^{(\eta)}$ in
order to avoid overloaded notation, when keeping dependencies on the occurring parameters.

\begin{theorem}\label{theo_dL_fin}
Let $\rho\in (0,\infty)$ and $d=d(L)$ such that $dL\to \theta \in [0,\infty)$.
If $\eta^{(L,N)}(0)$ is such that $\mu_{\#}\eta^{(L,N)}(0) \weakconv \mu_0\in E$, then
\begin{align}
    \left(\mu_{\#} \eta^{(L,N)}(t)\right)_{t\geq 0}\weakconv (\mu_t)_{t\geq 0}\,, \quad \text{in $D([0,\infty), E)$}\,,\text{ as }N/L\to \rho\,.
\end{align}
Here $(\mu_t)_{t\geq 0}$ denotes the measure-valued process on $E$
\eqref{eq_def_E} generated by $\genFin$, cf. \eqref{eq_def_L}, with initial
value $\mu_0$.
\end{theorem}

Note that the limit process does not depend on the density $\rho$ due to our choice of rescaling, which is discussed in more detail in Section~\ref{sec_overview}.  
In Proposition~\ref{prop_dLfin_pregenerator} we prove that the closure of $(\genFin,\mathcal{D}({\genFin}))$ is indeed the generator of a Feller process on the state space $E$.
\\

For the case $\theta=\infty$, i.e. $dL$ diverging, we expect clusters on the
scale $1/d$, cf. \cite{JCG19} and the discussion around \eqref{eq_exp_limi_JCG19} below.
Hence, in this case we consider the embedding
\begin{align}\label{eq_def_map_mv_dLinf}
    \Hat{\mu}^{(\cdot)}=\Hat{\mu}_{L,N}^{(\cdot)} :\Omega_{L,N} \to \mathcal{M}_1(\R_+)\quad\mbox{with}\quad \Hat{\mu}_{L,N}^{(\eta)}
    := \sum_{x=1}^L \frac{\eta_x}{N} \delta_{dL\frac{\eta_x}{N}} \,,
\end{align}
mapping particle configurations into the space of probability measures  $\mathcal{M}_1(\R_+)$, again with the topology induced by weak convergence of measures.
In contrast to $\theta<\infty$, any measure in $\mathcal{M}_1(\R_+)$ can be
approximated by particle configurations using \eqref{eq_def_map_mv_dLinf}, see
Lemma~\ref{lem_particle_config_approx_dLinf}, which is why we do not have to restrict ourselves to a strict subset of probability measures as above.

The lack of compactness of $\R_+$ now allows for diverging rescaled masses of
particle configurations in the thermodynamic limit, thus when $dL\to\infty$, we
expect the scaling limit to be a measure-valued process on
$\mathcal{M}_1(\overline{\R}_+)$, with $\overline{\R}_+=[0,\infty]$. We include
$\infty$ to describe mass on larger scales than $1/d$. 
Indeed the correct limit turns out to be a process on $\mathcal{M}_1(\overline{\R}_+)$ without interaction and single-particle operator
\begin{align}\label{eq_def_Ahat_single_part}
    {\mutInf}h(z) :=&
    (Bh)'(z)+
    (Bh(0)-Bh(z))\\
    =&
    z h''(z)+ (2-z)h'(z)+
    (h(0)-h(z))\nonumber
    \,,
\end{align}
acting on $h$ in the domain
\begin{align}\label{eq_def_adomain}
    \mathcal{D}({\mutInf}):= 
    \{
    h\, :\, h(\infty)=0 \text{ and } h|_{\R_+}\in C_c^3(\R_+)
    \}
    \cup \{ \text{constant functions} \}\subset C(\overline{\R}_+)\,.
\end{align}
We refer to the proof of Theorem~\ref{theo_theta_to_inf} for a derivation of $ \mutInf$ from $\mutFin$ in terms of a scaling
argument when $ \theta \to \infty$.

Slowing down the evolution of the inclusion process appropriately, we get the following result.

\begin{theorem}\label{theo_dL_inf}
Let $\rho\in (0,\infty)$ and $d=d(L)\to 0$ such that $dL\to \infty$. If
$\Hat{\mu}_{\#} \eta^{(L,N)}(0) \weakconv \Hat{\mu}_0\in \mathcal{M}_1(\overline{\R}_+)$, then 
\begin{align}\label{eq_theo_dLinf_conv}
   \left(\Hat{\mu}_{\#} \eta^{(L,N)} \big(\tfrac{t}{dL}\big)\right)_{t\geq 0} \weakconv (\Hat{\mu}_t)_{t\geq 0}\,,
   \quad \text{in $D([0,\infty),\mathcal{M}_1(\overline{\R}_+))$}
   \,, \text{ as }N/L\to \rho\,.
\end{align}
Here $(\Hat{\mu}_t)_{t\geq 0}$ denotes the measure-valued process on $\mathcal{M}(\overline{\R}_+)$ with initial value $\Hat{\mu}_0$ generated by
\begin{align}\label{eq_def_G}
    {\genInf}H(\mu)
    = 
    \sum_{1\leq k \leq n }
    \mu({\mutInf}h_k) \prod_{\substack{m=1 \\ m\neq k}}^n \mu(h_m)\quad\mbox{with}\  H(\mu)=\mu(h_1)\cdots \mu(h_n),\ h_k\in \mathcal{D}({\mutInf})\ .
\end{align}
\end{theorem}


The operator ${\genInf}$ 
may be interpreted as a Fleming-Viot process without interaction.
This is in contrast to the generator $\genFin$, which does not have
a Fleming-Viot interpretation, since the interaction term is not of the form 
\begin{equation*}
\begin{aligned}
2\sum_{1\leq k < l \leq n} \big( \mu(h_k h_l) -\mu(h_k) \mu(h_l)  \big)
    \prod_{j\neq k,l} \mu(h_j)\,,
\end{aligned}
\end{equation*}
see also Appendix~\ref{app_FV}.


We will see that the limiting dynamics are deterministic with absorbing state $\hat\mu =\mathrm{Exp}(1)\in \mathcal{M}_1(\R_+)$. 
In fact, the statement of Theorem~\ref{theo_dL_inf} can be reformulated into a hydrodynamic limit, cf. Proposition~\ref{prop_fokker_plank}.
Moreover, if $ \Hat{\mu}_0[{\R}_+]=1$, then $\Hat{\mu}_{t}[\R_+]=1$
for every $t\geq 0 $, cf. Corollary~\ref{cor_mass_process}, i.e. mass does not escape to larger scales.\\

The two theorems fully determine the dynamics of the inclusion process, with vanishing diffusivity, on complete graphs in the thermodynamic limit with density $\rho\in (0,\infty)$.
For a discussion of the boundary cases $\rho \in \{0,\infty\}$, we refer to
Section~\ref{sec_overview}. For measure valued processes with generators
\eqref{eq_def_L} and \eqref{eq_def_G} the evolution w.r.t. a simple test function (in the appropriate domain) is given by
\begin{equation}\label{eq_sdes}
d\mu_t (h)=\mu_t (\mutFin h)\, dt +dM_t^{(h)} \quad\mbox{and}\quad
d\Hat{\mu}_t (h)=\Hat{\mu}_t (\mutInf h)\, dt\,,
\end{equation}
respectively. 
Here $t\mapsto M_t^{(h)}$ is a martingale with (predictable) quadratic variation
\[
\big\langle M^{(h)}\big\rangle_t =\int_0^t \Big(\genFin\big(\mu_s (h)^2 \big) -2\mu_s (h)\genFin\mu_s (h)\Big)\, ds=\int_0^t \Big(\mu_s \big( (Bh)^2\big) -\mu_s
(Bh)^2\Big)\, ds\, ,
\]
given by the interaction term in $\genFin$ with $n=2$ and $h_1 =h_2=h$. Since $t\mapsto \hat{\mu}_t (h)$ solves a simple ODE without martingale part it is continuous for all $h\in\mathcal{D} (\hat A)$, so that the process $(\hat{\mu}_t)_{t\geq 0}$ has continuous paths. 
Continuity of the process $(\mu_t )_{t\geq 0}$, and thus of the martingale $(M_t^{(h)})_{t\geq 0}$, follows from the equivalence with the Poisson-Dirichlet diffusion, cf.\ Proposition~\ref{prop_equiv_pd_diff_superprocess}.

Taking expectations of the first term in \eqref{eq_sdes}, we see that $\Bar{\mu}_t :=\E_{\mu_0} [\mu_t ]$ satisfies
\[
\frac{d}{dt} \Bar{\mu}_t (h)=\Bar{\mu}_t (A_\theta h)\,, \quad\mbox{for all }h\in
C^2 ([0,1])\, ,
\]
which agrees with the time evolution of $(\mathbf{E}_{\mu_0} [h(Z(t)
)])_{t \geqslant 0}$ for a process $(Z(t) )_{t\geq 0}$ on $[0,1]$ with generator $A_\theta$ and initial distribution $\mu_0$. 
As a consequence, we have the following dualities
\begin{align}\label{eq_dual}
    \mathbb{E}_{\mu_0}[\mu_t (h)]=
    \mathbf{E}_{\mu_0}
    [h(Z(t))]
    \quad\mbox{and}\quad \Hat{\mu}_t (h) =\mathbf{E}_{\mu_0}
    [h(\hat{Z}(t))]\quad \forall t\geq 0\, ,
\end{align}
where $(\hat{Z}(t) )_{t\geq 0}$ is a process on $ \overline{\R}_{+}$ with
generator $\hat{A}$. In the latter case $\Hat{\mu}_t$ itself is deterministic
for fixed initial condition $\mu_0$ since it solves the ODE \eqref{eq_sdes}. Both processes $(Z(t) )_{t\geq 0}$ and $(\hat{Z}(t) )_{t\geq 0}$ are one-dimensional diffusions with resetting to $0$, which will be used in Sections \ref{sec_dLfin} and \ref{sec_dLinf} to study properties of the measure-valued processes.

\subsubsection{A size-biased viewpoint on the Poisson-Dirichlet diffusion} 

The process described in Theorem~\ref{theo_dL_fin} is a measure-valued process which provides an alternative description of the infinitely-many-neutral-alleles diffusion model introduced by Ethier and Kurtz in their seminal work \cite{EK81}.
Note that the process is also commonly referred to as Poisson-Dirichlet diffusion, which we will use throughout the paper.
The classical Poisson-Dirichlet diffusion is a Feller process on
$\overline{\nabla}$ with pre-generator\footnote{The original formulation of the pre-generator in \cite{EK81} includes a multiplicative factor of $\tfrac{1}{2}$ which we omitted here.}
\begin{align}\label{eq_def_PD_diffusion}
    {\genPD}f
    =
    \sum_{i,j=1}^\infty 
    p_i(\delta_{i,j}- p_j) \partial^2_{p_i p_j} f
    -\theta \sum_{i=1}^\infty p_i \partial_{p_i} f\,,
\end{align}
acting on functions in the domain
\begin{align}\label{eq_def_dmon}
    \mathcal{D}_{mon}({\genPD}):=
    \text{sub-algebra of } C(\overline{\nabla}) \text{ generated by }\ 
 1,\varphi_2,\varphi_3,\ldots \,,
\end{align}
where $\varphi_m (p):= \sum_{i=1}^\infty p_i^m$ for $m\geq 2$. ${\genPD}$ acts on such test functions with the convention that occurring sums on the r.h.s. of \eqref{eq_def_PD_diffusion} are evaluated on $\nabla := 
    \big\{ p\in\overline{\nabla}\,:\, \sum_{i=1}^\infty p_i= 1\big\}\subset \overline{\nabla}$ and extended to $\overline{\nabla}$ by continuity.
We stress that $\varphi_{1}(p):= \sum_{i =1}^{ \infty} p_{i}$ is not a
continuous function on $ \overline{\nabla}$.
Our size-biased approach circumvents such technical issues, which is one of its main advantages.

The name 
Poisson-Dirichlet diffusion is adequate, since its unique invariant distribution is the \emph{Poisson-Dirichlet distribution}
PD($\theta$).
The Poisson-Dirichlet distribution is a one-parameter family of probability measures supported on $\nabla$. 
It was first introduced by Kingman \cite{Ki75} as a natural limit of Dirichlet distributions. However, there is a more intuitive construction of the Poisson-Dirichlet distribution using a stick-breaking procedure, see for example \cite{Fe10}.
Later, the distribution was identified as the unique stationary measure of the split-merge dynamics \cite{DMWZZ04, Sc05} and the Poisson-Dirichlet diffusion \cite{EK81}.
Despite it being introduced in the field of population genetics, the Poisson-Dirichlet distribution has since then also appeared in statistical mechanics \cite{GUW11,BeUe11,IoTo20} and recently in interacting particle systems \cite{JCG19,CGG21}.

Naturally, one can consider the mapping of the Poisson-Dirichlet diffusion
under the isomorphism $\mu^{(\cdot)}$, cf. Lemma~\ref{lem_isomorphism_simplex_measures}, which yields a process on $E\subset\mathcal{M}_1([0,1])$.
Indeed, this push-forward process agrees with the process generated by
$\genFin$. The proof of this fact can be found in the appendix.

\begin{proposition}\label{prop_equiv_pd_diff_superprocess}
    Let $(\mu_t)_{t\geq 0}$ be the measure-valued process generated by
$\genFin$ \eqref{eq_def_L} with initial data $ \mu_{0} \in E $,
and 
$(X(t))_{t\geq 0}$ be the Poisson-Dirichlet
diffusion generated by $\genPD$ \eqref{eq_def_PD_diffusion} with $
\mu^{( X(0))} = \mu_{0}$, then 
    \begin{align}
       (\mu_t)_{t\geq 0}\weakeq \left(\mu^{(X(t))}\right)_{t\geq 0}  \,.
    \end{align} 
    In particular, the following properties translate immediately from $(X(t))_{t\geq 0}$ to $(\mu_t)_{t\geq 0}$:
    \begin{chosenEnum}
        \item The process $(\mu_t)_{t\geq 0}$ has a unique stationary distribution, which is reversible. It is given by $\mathbf{P}=\mu_\# \mathrm{PD}(\theta)$, i.e.   the law of
        \begin{align}
        \mu^{(X)}=
            \sum_{i=1}^\infty X_i \delta_{X_i}\,, \quad X\sim \mathrm{PD}(\theta)\,.
        \end{align}
        
        \item 
        The process $(\mu_t)_{t\geq 0}$ has continuous sample paths in $E$.

        \item For any initial value $\mu_0\in E$, we have
        \begin{align}
            \mathbb{P}(\mu_t(\{0\})=0 \ \forall t>0 ) = 1\,.
        \end{align}
   \end{chosenEnum}
\end{proposition}

Together with Theorem~\ref{theo_dL_fin}, this yields the following corollary. 

\begin{corollary}\label{cor_main_theo_dL_fin}
Let $\rho\in (0,\infty)$ and $d=d(L)$ such that $dL\to \theta\in[0,\infty)$. If $\eta^{(L,N)}(0)$ is such that $\tfrac{1}{N}\Hat{\eta}^{(L,N)}(0)\weakconv X(0)\in \overline{\nabla}$, then
\begin{align}
    \frac{1}{N}\left(\Hat{\eta}^{(L,N)}(t)\right)_{t\geq 0}\weakconv (X(t))_{t\geq 0}\,, \quad \text{ as }N/L\to \rho\,.
\end{align}
Here $(X(t))_{t\geq 0}$ denotes the Poisson-Dirichlet diffusion on
$\overline{\nabla}$ with parameter $\theta$ and initial value $X(0)$, generated
by ${\genPD}$, cf. \eqref{eq_def_PD_diffusion}, and $\hat\eta$ denotes the
ordered particle configuration.
\end{corollary}

\begin{remark}\label{rem_inv_prop1_3}
Proposition~\ref{prop_equiv_pd_diff_superprocess} yields a way to recover the
measure--valued process $( \mu_{t})_{t \geqslant 0}$ from the classical
Poisson-Dirichlet diffusion $( X_{t})_{t \geqslant 0}$. 
Since the map \eqref{eq_embedding_kingman_to_measure} is an isomorphism, see
Lemma~\ref{lem_isomorphism_simplex_measures}, and
has a continuous inverse ${(\mu^{(\cdot)})^{-1}:E \to \overline{\nabla}}$, 
the proposition also
implies $  ( X_{t})_{t \geqslant 0} \weakeq
\big( (\mu^{(\cdot)})^{-1}( \mu_{t}) \big)_{t \geqslant 0}$.
The inverse map
reconstructs elements in the Kingman simplex by counting and ordering mass occurrences, for
example, 
\begin{equation*}
\begin{aligned}
E \ni
\tfrac{1}{2}
\delta_{\frac{1}{4}} 
+ \tfrac{3}{8} \delta_{\frac{ 1}{8}}
+ \tfrac{1}{8} \delta_{0}
\mapsto 
\big( 
\tfrac{1}{4},
\tfrac{1}{4},
\tfrac{1}{8},
\tfrac{1}{8},
\tfrac{1}{8},
0, \ldots
\big)
\in \overline{\nabla}\,.
\end{aligned}
\end{equation*}
\end{remark}

In this article, we focus on a (joint) thermodynamic limit as $N,L\to\infty$ with $N/L\to\rho$,
but our derived scaling limits do not actually depend on the density $\rho$. Therefore our
approach also extends to different scaling regimes as is discussed in Section~\ref{sec_overview}.
The measure valued process generated by ${\genInf}$ is the natural extension of the process $\genFin$ (and thus to the Poisson-Dirichlet diffusion generated by ${\genPD}$) when $\theta \to \infty$. 
A first indication for this relationship can already be observed on the level of stationary distributions. 
From Proposition~\ref{prop_equiv_pd_diff_superprocess}(i), we recall that the stationary distribution w.r.t. $\genFin$ is given by the size-biased sample of $\mathrm{PD}(\theta)$. 
Consider $X^{(\theta)}\sim \mathrm{PD}(\theta)$ and sample an index $I \in \N$ such that
\begin{align}
    I=i\quad \text{with probability } X_i^{(\theta)}\,,
\end{align}
i.e. we pick the index $I$ with size-bias.
It is well known \cite[Theorem 2.7]{Fe10} that $X_I^{(\theta)}\sim \mathrm{Beta}(1,\theta)$. 
Moreover, Exp(1) is the absorbing state of the deterministic dynamics induced by ${\genInf}$, cf. Corollary~\ref{cor_exp_rever}. 
Now, the following connection between a Beta and an Exponential distribution holds: 
\begin{align}
        \theta\, \mathrm{Beta}(1,\theta) \weakconv \mathrm{Exp}(1)\,, \quad\mbox{as
}\theta\to\infty\,.
\end{align}
Hence, as $\theta\to\infty$, the rescaled size-biased sample $\theta\,
X^{(\theta)}_I$ converges weakly to an Exp(1) random variable.

This relationship can also be made sense of on the level of processes, summarised in the following diagram:
\[
\begin{tikzcd}
    && \big(\genFin,E\big) \arrow[r,Leftrightarrow]{..} \arrow{dd}{
    \theta \to \infty \ \text{when }\ \begin{cases} z\,
\mapsto\,  \theta z\,, \\  t\, \mapsto\,  t/\theta\,.\end{cases}} & \big({\genPD}, \overline{\nabla}\big)& \text{, if } \theta<\infty\\
    \text{ $\big(\mathfrak{L}_{L,N},\Omega_{L,N}\big)$} \arrow{r}{N/L\to \rho}&\Big\{ \arrow{ru}{\delta_{\frac{1}{N}\eta_x}} \arrow{rd}[swap]{\delta_{\frac{dL}{N}\eta_x}}&&&\\
    && \big({\genInf},\mathcal{M}_{ 1}(\overline{\R}_+)\big) &&\text{, if } \theta=\infty\,.
\end{tikzcd}
\]
We analyse the inclusion process $\big(\mathfrak{L}_{L,N},\Omega_{L,N}\big)$ and consider the two cases $\theta<\infty$ (Theorem~\ref{theo_dL_fin}) and $\theta=\infty$ (Theorem~\ref{theo_dL_inf}), with appropriate embeddings of configurations in the space of probability measures.
In the case $\theta<\infty$ the limiting process is equivalent to the Poisson-Dirichlet diffusion generated by ${\genPD}$ (Proposition~\ref{prop_equiv_pd_diff_superprocess}). 
Furthermore, our size-biased approach allows for a meaningful limit when $\theta=\infty$, identifying a natural extension for models with Poisson-Dirichlet diffusion limit, under appropriate rescaling of time and space. For this matter, we introduce the scaling operator $S_\theta: E \to \mathcal{M}_1(\R_+)$, which linearly scales measures on the unit interval to measures on the interval $[0,\theta]$, i.e. 
\begin{align}
    S_\theta: \mu(dz) \mapsto \mu(d\tfrac{z}{\theta})\,.
\end{align}

\begin{theorem}\label{theo_theta_to_inf}
Let $(\mu^\theta_t)_{t\geq 0}$ be the process generated by $\genFin$ and $(\Hat{\mu}_t)_{t\geq 0}$ be the process generated by ${\genInf}$. If $S_\theta \mu_0^\theta\weakconv \Hat{\mu}_0\in \mathcal{M}_1(\overline{\R}_+)$, then 
\begin{align}
    \big(S_\theta \mu_{t/\theta}^\theta  \big)_{t\geq 0}
    \weakconv (\Hat{\mu}_t)_{t\geq 0}\,,
    \quad \text{ in $C([0,\infty),\mathcal{M}_1(\overline{\R}_+))$}\,,
    \quad \text{ as $\theta \to \infty$}\,,
\end{align}
where we consider the topology induced by weak convergence on $\mathcal{M}_1(\overline{\R}_+)$.
\end{theorem}        
    
\begin{remark}
In the above theorem we saw that scaling of space (of order $\theta$) is
necessary to observe a meaningful limit as $\theta \to \infty$. Similarly, one
could scale the Kingman simplex to $\theta \overline{\nabla}$. However, in the
limit we lack the property of distinguishing between separate scales. Consider
for example $\theta = n^2$, $n\in \N$, and the sequence 
\begin{align}
    p^{(\theta)}
    = 
    \tfrac{1}{2}\big(
    \underbrace{\tfrac{1}{\sqrt{\theta}},\ldots, \tfrac{1}{\sqrt{\theta}}}_{\sqrt{\theta}\text{ times}},
    \underbrace{\tfrac{1}{\theta}, \ldots, \tfrac{1}{\theta}}_{\theta\text{ times}}, 0, \ldots
    \big)\in \nabla \,.
\end{align}
Then $\theta\, p^{(\theta)}_i\to\infty $ for every $i\in\N$.
On the other hand, 
\begin{align}
    S_\theta \mu^{(p^\theta)}
    = \tfrac{1}{2} \delta_{\frac{1}{2}} + \tfrac{1}{2}\delta_{\frac{1}{2}\sqrt{\theta}} \weakconv \tfrac{1}{2} \delta_{\frac{1}{2}} + \tfrac{1}{2}\delta_{\infty}\in \mathcal{M}_1(\overline{\R}_+)\,,
\end{align}
which captures both the amount of diverging mass and information on scales $\tfrac{1}{\theta}$. 
This highlights the fact that considering the space $E$, instead of
$\overline{\nabla}$, is essential for a detailed analysis of the
Poisson-Dirichlet diffusion in the boundary case $\theta \to \infty$.
In particular, this highlights that the gap in the schematic diagram
above (at the
bottom right) is not expected to be filled, by a naive scaling of $ \overline{\nabla}$.
\end{remark}

\vspace{.9cm}

\subsection{Comparison to the literature}

\subsubsection{Condensation and the inclusion process}

After its introduction \cite{GKR07}, the inclusion process has been subject to study as an interesting model of mass transport on its own \cite{CGGR13}.
In particular, in the context of condensation in stochastic particle systems it is a model of major interest. 
In short, a particle system exhibits condensation if 
a positive fraction of particles concentrates on sites with a vanishing volume fraction. 
Such sites with diverging occupation (and their occupying particles) are called
the \emph{condensate}, the remaining particles and sites are said to be the
\emph{background} or \emph{bulk} of the system.
In \cite{JCG19}, the existence of a non-trivial condensate was proven;
we refer to the same reference for an exact definition of the condensation phenomenon.
While the dynamics of the bulk for occupation numbers of order one is covered by general results on the propagation of chaos for particle systems \cite{GJ19}, 
the scope of the present article is to determine the dynamics of the mass distribution in the condensate. 
However, we want to stress that our analysis does not rely on  previous results on condensation. In fact, the clustering of particles on diverging scales is an implicit consequence of the scaling limits, presented in Theorems~\ref{theo_dL_fin} and~\ref{theo_dL_inf}.

For the inclusion process, the condensation phenomenon was first studied in \cite{GRV11}, however for spatially inhomogeneous systems on a finite lattice with diverging number of particles. 
In homogeneous systems, condensation is a consequence of increasing particle interactions relative to diffusion as $d\to 0$. 
The finite lattice case has been subject to further studies in
\cite{GRV13,BDG17,KS21, kim2021second}.
When considering a thermodynamic limit, i.e. diverging lattice size $L$ and
number of particles $N$ with finite limiting density $\rho$, condensation was
studied heuristically in \cite{CCG14}. They considered a one-dimensional
periodic lattice with totally asymmetric dynamics and vanishing diffusion rate
s.t. $\theta=0$. A modified model with stronger particle interactions, leading
to instantaneous condensation even in one spatial dimension has been considered
in \cite{WaEv12,CCG15}.

On a rigorous level, the thermodynamic limit of stationary distributions has been treated in \cite{JCG19}. Under the assumption that $d\to 0$ as $L\to\infty$, we have the following cases:
\begin{itemize}    
    \item if $dL\to 0$,  then the condensate is given by a single cluster, and if in addition $dL \log{L} \to 0$ this cluster is holding all the particles;


    \item if $dL\to \theta\in (0,\infty)$, then the condensate concentrates on macroscopic scales and is distributed according to a Poisson-Dirichlet distribution $\mathrm{PD}(\theta)$;

    \item on the other hand if $dL\to \infty$, the condensate is located on mesoscopic scales and the clusters are independent, more precisely,
    \begin{align}\label{eq_exp_limi_JCG19}
        d(\Tilde{\eta}_1, \ldots, \Tilde{\eta}_n) \weakconv \mathrm{Exp}(\tfrac{1}{\rho})^{\otimes n}\,.
    \end{align}
Here $\Tilde{\eta}$ denotes a size-biased sample w.r.t. $\eta\in \Omega_{L,N}$,
introduced in \eqref{eq_def_sb}.
\end{itemize}

The above result holds for any irreducible and spatially homogeneous dynamics on diverging finite graphs, where the inclusion process has stationary product measures.
Because condensation in homogeneous systems only occurs if $d\to 0$, \cite{JCG19} characterised the stationary condensates for the inclusion process, and in particular, the results hold for the complete graph dynamics we consider here.
It was proven in \cite{CGG21} that perturbations of the transition rates (in the case $\theta\leq 1$) still give rise to a Poisson-Dirichlet distributed condensate.

Both Theorem~\ref{theo_dL_fin} and Theorem~\ref{theo_dL_inf} complete the picture of condensation behaviour of the inclusion process on complete graphs outside of stationarity. Moreover, our results link the inclusion process dynamics directly with  the Poisson-Dirichlet diffusion, when $\theta<\infty$, which allows for enhanced understanding of the latter dynamics.

\subsubsection{The Poisson--Dirichlet diffusion}

In their work \cite{EK81}, Ethier and Kurtz proved that the closure of $({\genPD},\mathcal{D}_{mon}({\genPD}))$, cf. \eqref{eq_def_PD_diffusion}, gives rise to a generator of a diffusion process on $\overline{\nabla}$, which is the natural scaling limit of a finite-dimensional Wright-Fisher model when sending the number of individuals and types to infinity separately.
The restriction to test functions in $\mathcal{D}_{mon}({\genPD})$ turns out to be convenient, but makes it difficult to understand the precise dynamics of the infinite dimensional process. 
In \cite{EK81}, also an enlarged domain of test-functions of the form
\begin{align}
    p \mapsto \sum_{i=1}^\infty h(p_i)\,,\quad  h\in C^2([0,1]) \text{ with } h(0)=h'(0)=0\,,
\end{align}
was considered.
However, this does not improve the understanding of the dynamics on an intuitive level, which is particularly due to the fact that `sums are evaluated on $\nabla$ and extended to $\overline{\nabla}$ by continuity'.
In this paper, we instead propose to consider functions of the form 
\begin{align}\label{eq_def_h_mean_wrt_p}
    p\mapsto h(0) + \sum_{i=1}^\infty p_i (h(p_i)-h(0))\,,\quad h\in C^2([0,1])\,.
\end{align}
For a fixed $p\in \overline{\nabla}$, the r.h.s. is  the expectation w.r.t. the probability measure $\mu^{(p)}$, recall \eqref{eq_embedding_kingman_to_measure}.
Moreover, note that the functions $\varphi_m$ are of the form \eqref{eq_def_h_mean_wrt_p} with $h(p)=p^{m-1}$.

The usual approach in the literature, when constructing the Poisson-Dirichlet diffusion, is to take the large $L$-limit of an $L$-dimensional diffusion model. 
Alternatively, discrete models have been considered but then first convergence to the $L$-dimensional diffusion model, when $N\to \infty$, is proven. See for example \cite{EK81,CBERS17,RW09}. 
To the authors best knowledge, the present article is the first to consider a thermodynamic limit, which is taking both size of the system and number of particles to infinity at the same time while keeping the density approximately constant. 
This makes sense both from a physical and population genetics perspective.
In particular, taking a joint limit allows for interesting dynamics in the case $\theta = \infty$ which could otherwise not be considered, recall Theorem~\ref{theo_theta_to_inf}.

The Poisson-Dirichlet diffusion was treated previously as a measure-valued
process in $\mathcal{M}_1([0,1])$ in \cite{EG87,EK87}, where it was considered
as a Fleming-Viot process with mutation operator
    \begin{align}\label{eq:FV_mut}
        A_{FV}h(u) = \theta \int_0^1 [h(v)-h(u)]\,dv\,.
    \end{align}
Here the elements in $[0,1]$ are interpreted as types, and uniform jumps at rate $\theta$ in the mutation operator correspond to mutation events. These dynamics can be derived in the thermodynamic limit from the inclusion process on a complete graph with $dL\to\theta\in (0,\infty )$, using the embedding 
    \begin{align}\label{eq_FV_embedding}
        \eta \in \Omega_{L,N} \mapsto \sum_{x=1}^L \frac{\eta_x}{N} \delta_{\frac{x}{L}} \in\mathcal{M}([0,1])\,.
    \end{align}
This describes the spatial distribution of mass on the rescaled lattice and is
different from the approach in the present paper, where we ignore spatial
information and only keep track of the mass distribution, cf.\ \eqref{eq_embedding_kingman_to_measure}. On the other hand, our approach is more robust and allows for an extended analysis of the model for $\theta =\infty$. The embedding \eqref{eq_FV_embedding} has been considered in \cite{CT16} for the inclusion process on a complete graph of fixed size. They study convergence to equilibrium in the long-time limit and with diverging mass $N\to\infty$. 
The derivation of a Fleming-Viot process with mutation \eqref{eq:FV_mut} in the
thermodynamic limit is relatively straightforward if the inclusion process is
formulated in terms of particle positions, which is presented briefly in
Appendix~\ref{app_FV}.



Applying our approach to other geometries may be possible for dense random graphs along the lines of \cite{bovier2022}, which have diverging degrees leading to a self-averaging effect similar to the complete graph. 
In general, spatial models are difficult to treat since the inclusion process after the embedding \eqref{eq_embedding_kingman_to_measure} is not Markovian. Consider for example nearest-neighbour dynamics on a regular lattice, then it is known that the random-walk and the inclusion part of the dynamics have two different time scales, see \cite{ACR21}, and more sophisticated methods are necessary to treat this case.

\subsection{Outline of the paper}

In Section~\ref{sec_dLfin} we show that ${\genFin}$  generates a Feller process and prove Theorem~\ref{theo_dL_fin}.
We make use of explicit approximations of the inclusion process generator and the Trotter-Kurtz approximation theorem.
Moreover, we prove the equivalence of the Poisson-Dirichlet diffusion and our scaling limit in Section~\ref{sec_equiv}. 
Lastly, we discuss the advantages of considering size-biased dynamics in Section~\ref{sec_size_biased}. 
In Section~\ref{sec_dLinf} we determine the scaling limit when $\theta = \infty$, following a similar approach as in the case $\theta<\infty$.
We finish the section by proving the convergence $\tfrac{1}{\theta}\genFin\to {\genInf}$ stated in Theorem~\ref{theo_theta_to_inf}.
Lastly, we discuss boundary cases $\rho\in \{0,\infty\}$, fluctuations and open problems in Section~\ref{sec_overview}.

\section{Scaling limits in the case \texorpdfstring{$dL\to \theta<\infty$}{dL finite}}\label{sec_dLfin}
\subsection{The measure-valued process} 

In this section, we will prove that the measure-valued process generated by
$\genFin$ \eqref{eq_def_L} is a Feller process on the state space $E$
\eqref{eq_def_E}. Furthermore, we deduce weak convergence on the path space 
for  the inclusion process configurations embedded in the space of probability measures on the unit interval.

\subsubsection{Approximation of infinitesimal dynamics} 

The key result of this section is the following convergence result on the level of pre-generators
\begin{proposition}\label{prop_conv_gen_dL_fin}
Let $\rho\in (0,\infty)$ and $d=d(L)$ such that $dL\to \theta \in[0,\infty)$.
For every $H\in \mathcal{D}({\genFin})$, cf. \eqref{eq_def_DL}, we have with $\mathfrak{L}_{L,N}$ defined in \eqref{eq_def_ip_dynamics}
\begin{align}
    \lim_{N/L \to \rho} \sup_{\eta\in \Omega_{L,N}}
    \big|
    \mathfrak{L}_{L,N}H(\mu^{(\cdot)})(\eta)
    - (\genFin H)(\mu^{(\eta)})
      \big|=0\,.
\end{align}
\end{proposition}
We split the proof of Proposition~\ref{prop_conv_gen_dL_fin} into two parts. First, we only consider test functions of elementary form $H(\mu)=\mu(h)$, which corresponds to measuring a single observable $h\in C^3 ([0,1])$. We then extend the convergence result to arbitrary test functions in the domain, which requires to understand correlations between several observables. As usual, it turns out that only pairwise correlations contribute to leading order.

\begin{lemma}\label{lem_conv_single_part_dL_fin}
Let $\rho\in (0,\infty)$ and $d=d(L)$ such that $dL\to \theta \in[0,\infty)$.
Consider $H\in \mathcal{D}({\genFin})$ of the elementary form $H(\mu)=\mu(h)$, for some $h\in C^3([0,1])$. Then 
\begin{align}
    \lim_{N/L \to \rho} \sup_{\eta\in \Omega_{L,N}}
    \big|
    \mathfrak{L}_{L,N}H(\mu^{(\cdot)})(\eta)
    - \mu^{(\eta)}({\mutFin} h)
    \big|=0\,,
\end{align}
where ${\mutFin}$ is the single-particle generator, introduced in \eqref{eq_def_A_singlePart}.
\end{lemma}

\begin{proof}
Let $h\in C^3([0,1])$ and define $H(\mu):=\mu(h)$.
For the sake of convenience we introduce the notation $\Tilde{h}(z):= z\,
h(z)$. Thus, 
\begin{align}
H(\mu^{(\eta)})
=H(\mu_{\#} \eta) =
\mu^{(\eta)}(h)
=
\sum_{x=1}^L
\frac{\eta_x}{N}
{h}(\tfrac{\eta_x}{N})
=
\sum_{x=1}^L \Tilde{h}(\tfrac{\eta_x}{N})\,,
\end{align}
which allows us to write 
\begin{align}\label{eq_hTilde_Taylor_exp}
    H(\mu_{\#} \eta^{x,y})-H(\mu_{\#} \eta )
    &=
    \Tilde{h}(\tfrac{\eta_y+1}{N})- \Tilde{h}(\tfrac{\eta_y}{N})
+
\Tilde{h}(\tfrac{\eta_x-1}{N})-\Tilde{h}(\tfrac{\eta_x}{N}) \nonumber \\
&=\frac{1}{N}\Tilde{h}'(\tfrac{\eta_y}{N})
    +
    \frac{1}{2}\frac{1}{N^2}\Tilde{h}''(\tfrac{\eta_y}{N})
    -
    \frac{1}{N}\Tilde{h}'(\tfrac{\eta_x}{N})\\
    &\qquad +
    \frac{1}{2}\frac{1}{N^2}\Tilde{h}''(\tfrac{\eta_x}{N})
    +
    \frac{1}{6}
    \frac{1}{N^3} \Tilde{h}'''(\xi)\,,\nonumber
\end{align}
using a second-order Taylor approximation of $\Tilde{h}$, with $\xi\in [0,1]$.
Therefore, we have uniformly over configurations $\eta$
\begin{align*}
    \mathfrak{L}_{L,N}H(\mu^{(\cdot)})(\eta)
    =
    \sum_{\substack{x,y=1\\ x\neq y}}^{L}
    \eta_x (d+\eta_y)
    \Big[
    \frac{1}{N}\Tilde{h}'(\tfrac{\eta_y}{N})
    +
    &\frac{1}{2}\frac{1}{N^2}\Tilde{h}''(\tfrac{\eta_y}{N}) \\
    & -
    \frac{1}{N}\Tilde{h}'(\tfrac{\eta_x}{N})
    +
    \frac{1}{2}\frac{1}{N^2}\Tilde{h}''(\tfrac{\eta_x}{N})
    \Big]
    +o(1)\,,
\end{align*}
where $o(1)$ denotes a (uniformly in $ \Omega_{L,N}$) vanishing quantity  as $N/L
\to \rho$.
We split the sum into two parts, by analysing terms with coefficients $d\, \eta_x$ and $\eta_x \eta_y$ separately. We begin with the latter:
\begin{itemize}
    \item The contribution of inclusion rates $\eta_x \eta_y$ is limited to 
    \begin{align}\label{eq_lab_IP_2order}
        \sum_{\substack{x,y=1\\ x\neq y}}^{L}
        \frac{\eta_x}{N} \frac{\eta_y}{N}
        \Tilde{h}''(\tfrac{\eta_x}{N})
        = 
        \sum_{\substack{x=1}}^{L}\frac{\eta_x}{N} \Big(
        1-\frac{\eta_x}{N}
        \Big)
        \Tilde{h}''(\tfrac{\eta_x}{N})\,,
    \end{align}
    due to exact cancellation of the first-order terms $\Tilde{h}'$.
    
    \item On the other hand, contributions of the random-walk dynamics induced by rates $d \, \eta_x$ are given by
    \begin{align}\label{eq_lab_IP_1order}
        &d\, \sum_{\substack{x,y=1\\ x\neq y}}^{L}
        \frac{\eta_x}{N}
        \left[
        \widetilde{h}'(\tfrac{\eta_y}{N})
        +
        \frac{1}{2}\frac{1}{N}\widetilde{h}''(\tfrac{\eta_y}{N})
        -
        \widetilde{h}'(\tfrac{\eta_x}{N})
        +
        \frac{1}{2}\frac{1}{N}\widetilde{h}''(\tfrac{\eta_x}{N})
        \right]\nonumber \\
        &\qquad =
        dL\, \sum_{x=1}^{L} \frac{\eta_x}{N}
        \left[
        \frac{1}{L}
        \sum_{y\neq x}
        \widetilde{h}'(\tfrac{\eta_y}{N})
        -
        \frac{L-1}{L}
        \widetilde{h}'(\tfrac{\eta_x}{N})
        \right]
        +o(1)\,,
    \end{align}
    because second-order terms $\Tilde{h}''$ vanish in the thermodynamic limit due to 
    \begin{align}
        \left|
        d
        \sum_{\substack{x,y=1\\ x\neq y}}^{L}
        \frac{\eta_x}{2N^2}
        \big(\Tilde{h}''(\tfrac{\eta_x}{N})+\Tilde{h}''(\tfrac{\eta_y}{N}) \big)
        \right|
        \leq
         dL\,\sum_{x=1}^{L}\frac{\eta_x}{N^2} \|\Tilde{h}''\|_{\infty}
        \leq  \frac{dL}{N}  \|\Tilde{h}''\|_{\infty} \to 0\,.
    \end{align}
    Furthermore, we can absorb errors arising from replacing $\tfrac{L-1}{L}\Tilde{h}'$ with $\Tilde{h}'$, into $o(1)$. 
\end{itemize}
Now, combining \eqref{eq_lab_IP_2order} and \eqref{eq_lab_IP_1order} yields
\begin{align}\label{eq_disc_dyn_PD_diff_like_gen}
    \mathfrak{L}_{L,N}H(\mu^{(\cdot)})(\eta)
    =
    \sum_{\substack{x=1}}^{L}\frac{\eta_x}{N} \Big(
        1-\frac{\eta_x}{N}
        \Big)
        \Tilde{h}''(\tfrac{\eta_x}{N})
    +
    dL\, 
        \Big[
        \Tilde{h}'(0)
        -
        \sum_{x=1}^{L} \frac{\eta_x}{N}
        \Tilde{h}'(\tfrac{\eta_x}{N})
        \Big]
        +o(1)\,,
\end{align}
\sloppy where we additionally used Lemma~\ref{lem_riemann_partition_h0} to
conclude the uniform approximation $\frac{1}{L} \sum_{y=1\,, y\neq
x}^{L} \Tilde{h}'(\tfrac{\eta_y}{N})= \Tilde{h}'(0)+o(1)$.
Rewriting \eqref{eq_disc_dyn_PD_diff_like_gen} in terms of $\mu^{(\eta)}$, we have 
\begin{align}\label{eq_IP_disc_approx_A}
    \mathfrak{L}_{L,N}H(\mu^{(\cdot)})(\eta)
    =
    \mu^{(\eta)}
    \big(
    &z (1-z)h''(z) \nonumber \\
    &\qquad +2(1-z)h'(z)
    +dL (h(0)-h(z)-z h'(z))
    \big)+o(1)\,, 
\end{align}
where we used that $\Tilde{h}'(z)=h(z)+zh'(z)=Bh(z)$ and $\Tilde{h}''(z)=
2h'(z) +zh''(z)=(Bh)'(z)$.
Lastly, since $\|Bh\|_\infty<\infty$ and  $dL\to \theta$, we indeed have 
\begin{align}
    \mathfrak{L}_{L,N}H(\mu^{(\cdot)})(\eta)
    =\mu^{(\eta)}({\mutFin} h)+o(1)\,,
\end{align}
uniformly over all $\eta \in \Omega_{L,N}$. This concludes the proof.
\end{proof}

\begin{remark}
Note that the equivalence to the Poisson-Dirichlet diffusion can already be observed in \eqref{eq_disc_dyn_PD_diff_like_gen} when considering $h$ to be of the form $h(z)=z^{m-1}$, $m\geq 2$. In this case $H(\mu^{(\eta)})=\mu^{(\eta)}(h)=\sum_{x=1}^L \Tilde{h}(\tfrac{\eta_x}{N})=\varphi_m(\tfrac{\eta}{N})$, cf. \eqref{eq_def_dmon}, and
\begin{align}
     \mathfrak{L}_{L,N}H(\mu^{(\cdot)})(\eta)\simeq {\genPD}\varphi_m(\tfrac{\Hat{\eta}}{N})\,.
\end{align}
\end{remark}

After having proved the statement of Proposition~\ref{prop_conv_gen_dL_fin} for specific functions, we can now proceed with the proof of the full statement.

\begin{proof}[Proof of Proposition~\ref{prop_conv_gen_dL_fin}]
Let $H\in \mathcal{D}({\genFin})$. Without loss of generality we may assume $H$ has the form 
\begin{align}
    H(\mu) = \mu(h_1)\cdots \mu(h_n)\,, \quad h_k\in C^3([0,1])\,, \quad 1\leq k\leq n\,,
\end{align}
since linear combinations of such functions can be treated by linearity of the operators and the triangle inequality. Thus, considering $\eta\in \Omega_{L,N}$ and the configuration after  one particle jumped from $x$ to $y$, we have 
\begin{align*}
    H\left(\mu_{\#} \eta^{x,y}\right)
    &=
    \prod_{k=1}^n 
    \left(\mu_{\#} \eta^{x,y}\right)(h_k)\\
    &=
    \prod_{k=1}^n 
    \left[\widetilde{h}_k(\tfrac{\eta_y+1}{N})- \widetilde{h}_k(\tfrac{\eta_y}{N})
    +
    \widetilde{h}_k(\tfrac{\eta_x-1}{N})-\widetilde{h}_k(\tfrac{\eta_x}{N})+\mu^{(\eta)}(h_k)\right]\,.
\end{align*}
Now, expanding the product yields 
\begin{align}\label{eq_expansion_H_particle_shift}
    H\left(\mu_{\#} \eta^{x,y}\right)
    &=
    H\left(\mu_{\#} \eta\right)
    + \sum_{k=1}^n \left[\widetilde{h}_k(\tfrac{\eta_y+1}{N})- \widetilde{h}_k(\tfrac{\eta_y}{N})
    +
    \widetilde{h}_k(\tfrac{\eta_x-1}{N})-\widetilde{h}_k(\tfrac{\eta_x}{N})\right]\prod_{\substack{l=1\\l\neq k}}^n \mu^{(\eta)}(h_l)\nonumber \\
    &\qquad + 
    \sum_{1\leq k <l \leq n}
    \left[\widetilde{h}_k(\tfrac{\eta_y+1}{N})- \widetilde{h}_k(\tfrac{\eta_y}{N})
    +
    \widetilde{h}_k(\tfrac{\eta_x-1}{N})-\widetilde{h}_k(\tfrac{\eta_x}{N})\right]\\
    &\qquad \qquad \times 
   \left[\widetilde{h}_l(\tfrac{\eta_y+1}{N})- \widetilde{h}_l(\tfrac{\eta_y}{N})
    +
    \widetilde{h}_l(\tfrac{\eta_x-1}{N})-\widetilde{h}_l(\tfrac{\eta_x}{N})\right]
    \prod_{\substack{j=1\\j\neq k,l}}^n \mu^{(\eta)}(h_j)+ r(\eta)\,, \nonumber
\end{align}
with $r$ denoting the remainder. This expansion allows us to split
\begin{align}
    \mathfrak{L}_{L,N}H(\mu^{(\cdot)})(\eta)
    =
    \sum_{\substack{x,y=1\\x\neq y}}^L
\eta_x (d+\eta_y) \left[H(\mu_{\#} \eta^{x,y})-H(\mu_{\#} \eta )\right]
\end{align}
into three parts:
\begin{itemize}
    \item First, we make use of Lemma~\ref{lem_conv_single_part_dL_fin} which yields 
    \begin{align*}
        &\sum_{\substack{x,y=1\\x\neq y}}^L \eta_x (d+\eta_y)
        \sum_{k=1}^n [\widetilde{h}_k(\tfrac{\eta_y+1}{N})- \widetilde{h}_k(\tfrac{\eta_y}{N})
    +
    \widetilde{h}_k(\tfrac{\eta_x-1}{N})-\widetilde{h}_k(\tfrac{\eta_x}{N})]\prod_{\substack{l=1\\l\neq k}}^n \mu^{(\eta)}(h_l)\\
    &\qquad =
    \sum_{k=1}^n \mathfrak{L}_{L,N}\big(\mu^{(\cdot)}(h_k)\big)(\eta)  \prod_{\substack{l=1\\l\neq k}}^n \mu^{(\eta)}(h_l)
    =
    \sum_{k=1}^n \mu^{(\eta)}({\mutFin} h_k) \prod_{\substack{l=1\\l\neq k}}^n \mu^{(\eta)}(h_l)+o(1)\,.
    \end{align*}
    
    \item Next, we prove that the remainder $r$ has no contribution. More precisely, for any non-negative sequence $a_N$ satisfying $N^2 
    \, a_N\to 0$, i.e. $a_N$ lies in $o(\tfrac{1}{N^2})$, we have
    \begin{align}\label{eq_no_contr_of_small_terms}
        a_N
        \sum_{\substack{x,y=1\\x\neq y}}^L \eta_x (d+\eta_y) 
        \leq a_N \, (dL+N)N\to 0\,.
    \end{align}
    This includes, in particular, the remainder $r(\eta)$ because each summand lies in $o(\tfrac{1}{N^3})$, recall that each square bracket in \eqref{eq_expansion_H_particle_shift} vanishes uniformly like $N^{-1}$, cf. \eqref{eq_hTilde_Taylor_exp}.

    \item Lastly, we derive the interaction part where two observables are affected by the transition of a particle.
    Again, we perform a Taylor approximation for
    each of the two square brackets appearing in \eqref{eq_expansion_H_particle_shift}.
    Due to \eqref{eq_no_contr_of_small_terms}, together with \eqref{eq_hTilde_Taylor_exp}, it suffices to consider only products of  first-order terms $\Tilde{h}'$. Therefore, we are left with 
    \begin{align*}
        \frac{1}{N^2}
        \sum_{1\leq k <l \leq n}\sum_{\substack{x,y=1\\x\neq y}}^L \eta_x (d+\eta_y)
        \big[
     \Tilde{h}_k'(\tfrac{\eta_y}{N}) - \Tilde{h}_k'(\tfrac{\eta_x}{N})
     \big]\big[
     \Tilde{h}_l'(\tfrac{\eta_y}{N})- \Tilde{h}_l'(\tfrac{\eta_x}{N})
    \big]
        \prod_{\substack{j=1\\j\neq k,l}}^n \mu^{(\eta)}(h_j)+o(1)\,.
    \end{align*}
    For the same reason we include the random-walk interactions coming from
$d\, \eta_x$ in $o(1)$, and finally arrive at 
    \begin{align*}
        &\sum_{1\leq k <l \leq n}\sum_{\substack{x,y=1\\x\neq y}}^L \frac{\eta_x}{N} \frac{\eta_y}{N}
        \big[
     \Tilde{h}_k'(\tfrac{\eta_y}{N}) - \Tilde{h}_k'(\tfrac{\eta_x}{N})
     \big]\big[
     \Tilde{h}_l'(\tfrac{\eta_y}{N})- \Tilde{h}_l'(\tfrac{\eta_x}{N})
    \big]
        \prod_{\substack{j=1\\j\neq k,l}}^n \mu^{(\eta)}(h_j)+o(1)\\
        &\qquad =
        2\sum_{1\leq k <l \leq n}
    \big(
    \mu^{(\eta)}(\Tilde{h}'_k \Tilde{h}'_l)
    -
    \mu^{(\eta)}(\Tilde{h}'_k)\mu^{(\eta)}(\Tilde{h}'_l)
    \big)
        \prod_{\substack{j=1\\j\neq k,l}}^n \mu^{(\eta)}(h_j)+o(1)\,,
    \end{align*}
    where we expanded the product of square brackets and added the (non-contributing) diagonal $x=y$, before writing the expression in terms of $\mu^{(\eta)}$. Also, recall that $\Tilde{h}'=Bh$.
\end{itemize}

Overall, combining the three bullets above, we derive
\begin{align}
    \mathfrak{L}_{L,N}H(\mu^{(\cdot)})(\eta)
    =
    \genFin H(\mu^{(\eta)}) +o(1)\,, 
\end{align}
uniformly in $\eta\in \Omega_{L,N}$. This finishes the proof.
\end{proof}

\subsubsection{Convergence to the measure-valued process}

The measure-valued process takes values in the space $E= 
    \mu^{(\overline{\nabla})}\subset \mathcal{M}_1([0,1])$, cf. \eqref{eq_def_E}. 
Due to Lemma~\ref{lem_isomorphism_simplex_measures} $E$ itself is closed, thus, compact w.r.t. the topology induced by weak convergence of probability measures.
Because this topology coincides with the subspace topology, the Hausdorff property of $E$ is inherited from $\mathcal{M}_1([0,1])$.

In this section, we show that the dynamics described by ${\genFin}$ give rise to a Feller process and prove Theorem~\ref{theo_dL_fin}, which states that the process arises naturally as the scaling limit of the inclusion process.

\begin{proposition}\label{prop_dLfin_pregenerator}
For $\theta\in [0,\infty)$
the linear operator $(\genFin,\mathcal{D}({\genFin}))$ is closable and its closure generates a Feller process on the state space $E\subset \mathcal{M}_1([0,1])$.
\end{proposition}

The proof follows along the lines of \cite[Theorem 2.5]{EK81} where they proved existence of the Poisson-Dirichlet diffusion. 

\begin{proof}
Throughout the proof we will make use of the  sub-domain
\begin{align}
    \mathcal{D}_{mon}({\genFin}):= & \Big\{
    \text{sub-algebra of } C(E) \text{ generated by functions}\nonumber\\ &\qquad \text{ $\mu \mapsto
\mu(h)$  with }h(z)=z^m\,, m\in \N_0 \Big\}
    \subset  \mathcal{D}({\genFin})\,
\label{eq_def_dmonL}.
\end{align}
First note that, due to the Stone-Weierstrass theorem, $\mathcal{D}_{mon}({\genFin})$ (and therefore $\mathcal{D}({\genFin})$) is dense in $C(E)$ since it separates points: consider $\mu,\sigma\in E$ such that $\mu\neq \sigma$, then $\mu(z^m)\neq \sigma(z^m)$ for some $m\in \N$ since otherwise all moments, and hence $\mu$ and $\sigma$, agree.

Next, dissipativity of $(\genFin ,\mathcal{D}({\genFin}))$ follows from that of $(\mathfrak{L}_{L,N})_{L,N}$, since for any $H\in \mathcal{D}({\genFin})$ we have 
\begin{align}
    \|(\lambda-
\mathfrak{L}_{L,N})H(\mu^{(\cdot)})\|_{\Omega_{L,N},\infty} \geq \lambda \|H(\mu^{(\cdot)})\|_{\Omega_{L,N},\infty}\,\quad \forall \lambda>0.
\end{align}
The left hand side is upper bounded by 
\begin{align}
    \|(\lambda-
\genFin )H\|_{E,\infty}
+
\|\genFin H(\mu^{(\cdot)})-
\mathfrak{L}_{L,N}H(\mu^{(\cdot)})\|_{\Omega_{L,N},\infty}\,,
\end{align}
with the second term vanishing due to Proposition~\ref{prop_conv_gen_dL_fin}. On the other hand, using Lemma~\ref{lem_particle_config_approx_simplex}, we have 
\begin{align}
    \sup_{\eta\in \Omega_{L,N}}|H(\mu^{(\eta)})|\to \sup_{p\in \overline{\nabla}} |H(\mu^{(p)})| = \|H\|_{E,\infty}\,.
\end{align}

In the remainder of the proof, we first conclude that $\mathcal{D}_{mon}({\genFin})$ is a core for $\genFin$, using the fact that $\genFin$ is triangulisable.
The full statement then follows immediately by an extension argument.
For that purpose, we define subspaces
\begin{align}
    D_n({\genFin}):=
    \{
    H\in \mathcal{D}_{mon}({\genFin}) \,: \, deg(H)\leq n
    \},
\end{align}
where $deg(H)=m_1+\cdots + m_k$ if $H$ is of the form $\mu(z^{m_1})\cdots\mu(z^{m_k})$, $m_j\in\N$ for $1\leq j\leq k$. When $H$ is given by linear combinations of such products, the degree denotes the maximum degree of the products. 
Note that $\big( D_n({\genFin})\big)_{n\geq 1}$ defines an increasing sequence with limit $\mathcal{D}_{mon}({\genFin})$. 
It is only left to show that $\genFin$ maps elements of $D_n({\genFin})$ back into itself. This is, however, immediate since both parts of the generator $\genFin$ \eqref{eq_def_L}
map polynomials of a certain degree back into polynomials of the same degree.
Hence, using \cite[Proposition I.3.5]{EK_book}, we conclude that $(\genFin ,\mathcal{D}_{mon}({\genFin}))$ is indeed closable and gives rise to a strongly continuous contraction semigroup $(T_t)_{t\geq 0}$ on $C(E)$. 

Now, we can easily verify that also $\mathcal{D}({\genFin})$ is a core, using
\cite[Proposition I.3.1]{EK_book}, since 
\begin{align}
    \mathcal{R}(\lambda-{\genFin}|_{\mathcal{D}({\genFin})})
    \supset
    \mathcal{R}(\lambda-{\genFin}|_{\mathcal{D}_{mon}({\genFin})})\,,
\end{align}
is dense for some $\lambda>0$. Since generators are maximal dissipative, we
know that the closures w.r.t. both cores must agree, cf. \cite[Proposition
I.4.1]{EK_book} and hence  give rise to the same semigroup $(T_t)_{t\geq 0}$.
It is only left to show that the semigroup is positive and conservative, in
particular $E$ is invariant under the dynamics $\genFin$. In order to see this,
we apply Trotter's theorem, see e.g. \cite[Theorem I.6.1]{EK_book}, which concludes that Proposition~\ref{prop_conv_gen_dL_fin} implies 
\begin{align}\label{eq_trotter}
\lim_{N/L\to\rho}
    \sup_{\eta \in \Omega_{L,N}}
    \big|
    \mathfrak{T}^{(L,N)}_t H(\mu^{(\cdot)})(\eta)
    - (T_t H)(\mu^{(\eta)})
    \big|=0
    \,,\quad \forall H\in C(E)\,, t\geq 0\,,
\end{align}
where $\mathfrak{T}^{(L,N)}$ denotes the semigroup generated by
$\mathfrak{L}_{L,N}$. Now, both positivity and conservation follow from those
of $(\mathfrak{T}^{(L,N)})_{L,N}$. This conludes the proof.
\end{proof}

\begin{remark}
It is natural to ask why one should go through the inconveniences of extending the core from $\mathcal{D}_{mon}({\genFin})$ to $\mathcal{D}({\genFin})$. However, we will see in the next subsection that the extended core allows for a better interpretation of the underlying dynamics in the Poisson-Dirichlet diffusion.
\end{remark}

The proof of our first main result, namely, convergence of the inclusion process (when embedded in the space of probability measures) to the measure-valued process characterised by $\genFin$, is now an immediate consequence of a classical convergence theorem.

\begin{proof}[Proof of Theorem~\ref{theo_dL_fin}]
We apply \cite[Theorem IV.2.11]{EK_book} together with \eqref{eq_trotter}, which immediately concludes the desired convergence result. 
\end{proof}

Finally, we can use Theorem~\ref{theo_dL_fin} to prove convergence of the inclusion process to the Poisson-Dirichlet diffusion. 

\begin{proof}[Proof of Corollary~\ref{cor_main_theo_dL_fin}]
Every function $\varphi_m$ can be written in terms of an expectation $$p\mapsto
\mu^{(p)}(h_m)=\varphi_{m}(p)\,,$$
with $h_m(z):=z^{m-1}$. Thus, we have
\begin{align}\label{eq_supp_PD_conv}
    \Big(\varphi_m\big(\tfrac{1}{N}\Hat{\eta}^{(L,N)}(t)\big)\Big)_{t\geq 0}
    = \Big( (\mu_{\#} \eta^{(L,N)} (t))(h_m) \Big)_{t\geq 0}
    \weakconv
    \big(\mu_t(h_m)\big)_{t\geq 0}\,,
\end{align}
using Theorem~\ref{theo_dL_fin}.
This convergence can be extended to arbitrary elements in
$\mathcal{D}_{mon}({\genPD})$, in particular such sequences are tight in $D([0,
\infty), \R)$. Thus, 
 the sequence
$\big(\big(\tfrac{1}{N}\Hat{\eta}^{(L,N)}(t) \big)_{t\geq 0}\big)_{L,N}$ is
tight in $D([0,\infty), \overline{\nabla})$ and has subsequential limits, see
e.g. \cite[Theorem III.9.1]{EK_book}.
As convergence of finite dimensional marginals follows from \eqref{eq_supp_PD_conv}, we conclude the statement together with Proposition \ref{prop_equiv_pd_diff_superprocess}, which is proved in the next subsection.
\end{proof}

\subsection{Equivalence of the measure-valued process with PD-diffusion}\label{sec_equiv}

In this section we prove Proposition \ref{prop_equiv_pd_diff_superprocess} and investigate the equivalence of the measure-valued process generated by $\genFin$ \eqref{eq_def_L} and the Poisson-Dirichlet diffusion on the simplex $\overline{\nabla}$, generated by ${\genPD}$ \eqref{eq_def_PD_diffusion}.
We already saw in the proof of Lemma~\ref{lem_conv_single_part_dL_fin}, cf. \eqref{eq_disc_dyn_PD_diff_like_gen}, the similarity of dynamics of $\genFin$ and ${\genPD}$. Indeed, a simple calculation shows that the two can be linked: Using the embedding \eqref{eq_embedding_kingman_to_measure} we get for all 
$p\in \overline{\nabla}$ and $H(\mu)=\mu(h)$, with $h\in \mathcal{D}(A)$,
\begin{align*}
    \genFin H(\mu^{(p)})
    = &
    \mu^{(p)}(\mutFin h)
    =  \big(
    1-\| p\|_1 
    \big)\mutFin h(0) \\
    & + \sum_{i=1}^\infty p_i \Big(
    p_i(1-p_i) h''(p_i)+\big( 2(1-p_i)-\theta p_i \big) h'(p_i) + \theta(h(0)-h(p_i)) \Big)\,.
\end{align*}
Defining now $f(p):= \mu^{(p)}(h)$, we have 
\begin{align}\label{eq_EK_PD_A+_operator}
    \genFin H(\mu^{(p)})
    =
    2h'(0)(
    1-\|p\|_1
    )
    +
    {\genPD} f(p)\,,
\end{align}
where we used that $\partial_{p_i} f(p) = -h(0) +p_i h'(p_i)+h(p_i ))$ and $\partial_{p_i p_j} f(p)
= \delta_{ij} \big( 2h'(p_i) +p_i h''(p_i)\big)$ and $\genPD$ as defined in \eqref{eq_def_PD_diffusion}.

\begin{remark}
In \cite{EK81}, the authors extended the domain of ${\genPD}$ from $\mathcal{D}({\genPD})$ to the sub-algebra of $C(\overline{\nabla})$ generated by functions of the form $p\mapsto \sum_{i=1}^\infty g(p_i)$, with $g\in C^2([0,1])$ such that $g(0)=g'(0)=0$. 
This yields a similar expression as \eqref{eq_EK_PD_A+_operator}, cf. \cite[Display (2.17)]{EK81}. However, the expression again only made sense with the convention that sums are evaluated on $\nabla$ and extended by continuity, in which case the first summand in \eqref{eq_EK_PD_A+_operator} disappears.
\end{remark}

\begin{proof}[Proof of Proposition~\ref{prop_equiv_pd_diff_superprocess}]
In order to show the equivalence of the two processes, it suffices to restrict ourselves to the domains generated by monomials as defined in \eqref{eq_def_dmon} and \eqref{eq_def_dmonL}. 
Every function $H\in \mathcal{D}_{mon}({\genFin})$ can be mapped to $f_H\in \mathcal{D}_{mon}({\genPD})$ (and vise versa).
Let $H\in \mathcal{D}_{mon}({\genFin})$ be of the form $H(\mu) = \mu(h_1)\cdots \mu(h_n)$, where $h_k(z):=z^{m_k-1}$, then  $f_H=\varphi_{m_1}\cdots \varphi_{m_n}$ where we recall $\varphi_{m} (p)=\sum_{i=1}^\infty p_i^m$.
Then 
\begin{align}
    ({\genFin} H)(\mu^{(p)})
    &=
    2\sum_{1\leq k < l \leq n} m_k m_l \big( \mu^{(p)}(h_k h_l) -\mu^{(p)}(h_k) \mu^{(p)}(h_l)  \big)
    \prod_{j\neq k,l} \mu^{(p)}(h_j)\\
    &\qquad +\sum_{1\leq k \leq n} \mu^{(p)}(\mutFin h_k)\prod_{j\neq k}
\mu^{(p)}(h_j)\,, \nonumber
\end{align}
where we used that $Bh_k=m_k\, h_k$. Rewriting the r.h.s. in terms of $\varphi$'s, we have 
\begin{align*}
    ({\genFin} H)(\mu^{(\cdot)})
    &=
    2\sum_{1\leq k < l \leq n} m_k m_l \big( \varphi_{m_k+m_l-1} -\varphi_{m_k} \varphi_{m_l} \big)
    \prod_{j\neq k,l} \varphi_{m_j}\\ &\qquad +\sum_{1\leq k \leq n} {\genPD}\varphi_{m_k} \prod_{j\neq k} \varphi_{m_j}\,,
\end{align*}
where we used $\mu^{(p)}(\mutFin h_k)={\genPD} \varphi_{m_k}$ from
\eqref{eq_EK_PD_A+_operator} and the fact that $\varphi_{m_{k}} \in \mathcal{D}_{mon}({\genPD}) $ since
$m_{k} \neq 1 $.
Thus, $({\genFin} H)(\mu^{(\cdot)})$ agrees with ${\genPD} f_H$ on $\mathcal{D}_{mon}({\genPD})$, cf. \cite[Display (2.13)]{EK81}.
Let $(X(t))_{t\geq 0}$ be the Poisson-Dirichlet diffusion, then for every $H\in \mathcal{D}_{mon}({\genFin})$ 
\begin{align}
    H(\mu^{(X(t))}) -\int_0^t {\genFin} H(\mu^{(X_s)}) \, ds
    &= f_H(X(t))-  \int_0^t {\genPD} f_H(X_s) \, ds
\end{align}
defines a martingale in $t$. Thus, $(\mu^{(X(t))})_{t\geq 0}$ solves the martingale problem for $({\genFin}, \mathcal{D}_{mon}({\genFin}))$.

Now, it is almost immediate that properties (i) -- (iii) in Proposition~\ref{prop_equiv_pd_diff_superprocess} hold for the measure valued process $(\mu_t )_{t\geq 0}$. 
First, let $G,H\in \mathcal{D}_{mon}({\genFin})$ and choose  corresponding $f_G, f_H\in \mathcal{D}_{mon}({\genPD})$ as above. 
We know that $({\genFin}G)(\mu^{(p)})={\genPD}f_G(p)$. 
Writing $\nu = \mathrm{PD}(\theta)$ for simplicity, we have $\mathbf{P}=\mu_{\#} \mathrm{PD}(\theta )$ with
\begin{align}
    \mathbf{P}(H{\genFin}G)
    =
    \nu \big(H(\mu^{(p)})({\genFin}G)(\mu^{(p)})\big)
    =
    \nu (f_H{\genPD}f_G)\,.
\end{align}
It is known that $\mathrm{PD}(\theta)$ is the unique invariant distribution w.r.t. ${\genPD}$ \cite[Theorem 4.3]{EK81}; which is also reversible. Together with the above display, this yields $\mathbf{P}(H{\genFin}G)=\mathbf{P}(G{\genFin}H)$.

Continuity of the trajectories in (ii) follows from the diffusion property of $(X(t)(\omega))_{t\geq 0}$ and continuity of the map $\mu^{(\cdot)}$, together with the fact $(\mu^{(X(t))})_{t\geq 0} \stackrel{d}{=} (\mu_t)_{t\geq 0}$.

Lastly, (iii) is a consequence of
\begin{align*}
    \mathbb{P}\left[\mu_t(\{0\})=0 \quad \forall t>0 \right] 
    =
    \mathbb{P}\left[\mu^{(X(t))}(\{0\})=0 \quad \forall t>0 \right]
    =
    \mathbb{P}\left[X(t) \in 
   \nabla \quad \forall t>0 \right] =1\,,
\end{align*}
where we used \cite[Theorem 2.6]{EK81} in the last step.
\end{proof}

\begin{remark}
Naturally, we could have proven convergence of the inclusion process to the Poisson-Dirichlet diffusion directly and then defined the measure-valued dynamics using the embedding via $\mu^{(\cdot)}$. This would have slightly shortened the exposition in the present section, since it would not have been necessary to verify existence of the limiting dynamics.
We refrained from doing so for the sake of a better understanding of the underlying dynamics in the measure-valued process, in particular on the extended domain $\mathcal{D}({\genFin})$.
\end{remark}

Recall the generator of the single-particle dynamic \eqref{eq_def_A_singlePart}
\begin{align*}
    {\mutFin} h(z):= z(1-z) h''(z) +(2(1-z)-\theta z)h'(z) +\theta(h(0)-h(z))\,, \quad h\in C^2([0,1])\,,
\end{align*}
which characterises a Feller process on the unit interval. The process evolves according to a diffusion with an additional renewal mechanism due to jumps to zero.

\begin{lemma}
    The Beta distribution $\mathrm{Beta}(1,\theta)$ is the unique invariant distribution with respect to  ${\mutFin}$. 
\end{lemma}

\begin{proof}
    For $\theta=0$, we interpret the degenerate distribution $\mathrm{Beta}(1,0)$ as the Dirac point mass $\delta_1$.
    The statement is then clear since $\mutFin h(1)=0$. Hence, in this case with $\theta=0$ the point mass is even reversible.
    
    Now, let $\theta>0$ and consider $H(\mu):=\mu(h)$, then by Proposition~\ref{prop_equiv_pd_diff_superprocess}(i)
    \begin{align}
        0=\mathbf{P}({\genFin} H)
        = \int  \mu(\mutFin h)\, \mathbf{P}(d\mu)
        =
        \E\Big[
        \sum_{i=1}^\infty X_i \mutFin h(X_i)
        \Big]
        =
        \E[\mutFin h(\Tilde{X}_1)]\,,
    \end{align}
    where $X\sim \mathrm{PD}(\theta)$. It is well known that the first size-biased marginal $\Tilde{X}_1$ is $\mathrm{Beta}(1,\theta)$-distributed \cite[Theorem 2.7]{Fe10}.
    
    Uniqueness of the Beta distribution is due to Harris recurrence of the process, see e.g. \cite{MeTw93}.
    For the case $\theta>0$, the resetting mechanism guarantees that the process returns to zero infinitely often  almost surely.
    On the other hand for $\theta=0$, $\mutFin$ agrees with a Jacobi diffusion, see \eqref{eq_def_jacobi_diff} below. The corresponding process runs into the absorbing state $z=1$ in finite time, independent of the initial condition.
\end{proof}

Due to the jumps to zero, one does not expect that $\mathrm{Beta}(1,\theta)$, $\theta>0$, is reversible w.r.t. ${\mutFin}$. Indeed, this can be  verified easily by considering the example $h(x)=x$ and $g(x)=x^2$, in which case   
\begin{align}
    \mathrm{Beta}(1,\theta)(g \mutFin h)
    =
    \frac{8\theta \Gamma (\theta+1)}{\Gamma(\theta+4)}
    \neq 
    \frac{6\theta \Gamma (\theta+1)}{\Gamma(\theta+4)}
    =
    \mathrm{Beta}(1,\theta)(h \mutFin g)\,, \quad \forall \theta>0\,.
\end{align}

\subsection{The advantage of a size-biased evolution}\label{sec_size_biased}

Throughout the previous sections, we have seen two viewpoints of the same dynamics. The classical Poisson-Dirichlet diffusion considers a ranked configuration space.
However, this obscures the dynamics on microscopic scales, which results e.g.
into defining the r.h.s. of the generator $\genPD$ \eqref{eq_def_PD_diffusion}
to be \emph{evaluated on $\nabla$ and extended to $\overline{\nabla}$ by
continuity}, because $\varphi_{1}= \| \cdot \|_{1}$ does not lie in the domain
$\mathcal{D}_{mon}({\genPD})$. 
Alternatively, one can consider unordered dynamics, i.e. observing the evolution from a fixed position, or with a size-biased viewpoint. 
The Poisson-Dirichlet diffusion concentrates immediately on configurations consisting of macroscopic-sized fragments, which can only concentrate on a vanishing fraction of the volume. Hence, an unordered state space can only describe dynamics up to a certain point when the mass present at the observed positions disappears.
The goal of this section is to emphasise that a size-biased viewpoint allows for both, a complete description of the macroscopic dynamics while observing interaction with the microscopic scale.


\subsubsection{Time evolution on fixed sites}\label{sec_cf_sb_fix_site}

First, we look at arbitrary finite positions and observe the evolution of masses on them. As the inclusion process is spatially homogeneous, we may choose for simplicity $\eta \mapsto (\eta_1,\ldots, \eta_n)$.

We start by only considering the evolution on the first site. Performing similar approximations as in Section~\ref{sec_dLfin}, we can see for an arbitrary function $h\in C^{3}([0,1])$
\begin{align}
    \mathfrak{L}_{L,N}h(\tfrac{(\cdot)_1}{N})(\eta)
    =
    A_{Jac(\theta)}h(\tfrac{\eta_1}{N})+o(1)\,,
\end{align}
where 
\begin{align}\label{eq_def_jacobi_diff}
    A_{Jac(\theta)}h(z):= z(1-z)h''(z) -\theta z\, h'(z)\,.
\end{align}
In fact, we can see $A_{Jac(\theta)}$ emerging in \eqref{eq_disc_dyn_PD_diff_like_gen} when fixing a position $x$.
The operator $A_{Jac(\theta)}$ is the generator of a Jacobi-diffusion, cf. \cite{FSRW21},  and describes the evolution of a single chunk of mass located at a given position.

To describe the evolution on the first $n$ positions
we introduce for $i=0,1,\ldots ,n$
\begin{align}
    \xi_i = \xi_i (\eta) :=
    \begin{cases}
    \eta_i \quad & \text{if}\; 1\leq i\leq n\\
    N-\sum_{j=1}^n{\xi_j} \quad & \text{if}\; i=0\,,
    \end{cases}
\end{align}
where $\xi_0$ is the remaining mass in the system outside sites $1,\ldots ,n$. 
Thus, the rescaled vector $\tfrac{1}{N}\xi$ lies in $\Delta_{n+1}:=\{p \in [0,1]^{n+1}\,:\, \sum_{i=0}^n p_i = 1\}$.
Again, by approximation of the generators, one can show that
\begin{align}
    \tfrac{1}{N}(\xi(t))_{t\geq 0} \weakconv 
    \mathrm{WF}_{n+1}(\theta, 0, \ldots, 0)\,.
\end{align}
Here, $\mathrm{WF}_{n+1}(\theta, 0, \ldots, 0)$ denotes the Wright-Fisher diffusion on $\Delta_{n+1}$ which is characterised by the generator 
\begin{align}\label{eq_WF}
    A_{WF_{n+1}(\theta,\mathbf{0})}
    h(z_{0,\dots, n})
    &=
    \sum_{i,j=0}^n
    z_i (\delta_{i,j}- z_j) \partial^2_{z_i z_j }h(z_{0,\dots, n})\\
    &\qquad + \theta \sum_{i=1}^n z_i (\partial_{z_0} h-\partial_{z_i} h)(z_{0,\dots, n})\,, \nonumber 
\end{align}
acting on those $h$ that have an extension to $\R^{n+1}$ which is twice continuously differentiable.

However, we can already see for a single observable that the Jacobi-diffusion has an absorbing state at $z=0$, which it will run into in finite time almost surely \cite[Section 7.10]{Dur08}.
Similarly, the Wright-Fisher diffusion will be absorbed at $(1,0, \ldots 0)\in \Delta_{n+1}$, after which the process does not capture the dynamics of the infinite-dimensional process anymore as all the mass has moved away from the first $n$ sites.

For the Poisson-Dirichlet diffusion, the relationship to the Jacobi and Wright-Fisher diffusion has been studied in greater generality in the two-parameter setting \cite{FSRW21}. 
They use a Fleming-Viot construction of the process, cf. \eqref{eq_FV_embedding}. Because they start from the Poisson-Dirichlet diffusion on $\overline{\nabla}$, there is no underlying graph structure and instead of placing mass at a fixed position, they choose a uniform random variable on $[0,1]$ which determines the position of the point mass. 

Moreover, it is interesting to note that the boundary behaviour of the Jacobi diffusion agrees with the one of the PD-diffusion. More precisely, in the case of the Jacobi diffusion the state $1$ can only be reached if and only if $\theta<1$, see e.g. \cite[Theorem 4.1]{Shi81} or \cite[Section 7.10]{Dur08}. Similarly, the PD-diffusion $(X(t))_{t\geq 0}$ hits the finite dimensional sub-simplices $\nabla \cap \{\sum_{i=1}^n p_i=1\}$, for any $n\geq 1$, if and only if $\theta<1$ \cite{Sc91}.



\subsubsection{Duality and size-biased time evolution}

We recall from Dynkin's formula, i.e. taking expectations of the first term in \eqref{eq_sdes},
\[
\frac{d}{dt} \E_{\mu_0}\big[ \mu_t (h)\big] =\E_{\mu_0}\big[ \mu_t (A_\theta h)\big]\,, \quad\forall h\in
 \mathcal{D}(A)\, ,
\]
which implies the duality \eqref{eq_dual}
\begin{align}\label{eq_dual2}
    \mathbb{E}_{\mu_0}[\mu_t (h)]=
    \mathbf{E}_{\mu_0}
    [h(Z(t))]\,, \quad \forall h \in \mathcal{D}(A)\,,
\end{align}
where $(Z(t) )_{t\geq 0}$ is a process on $[0,1]$ with generator $A_\theta$, cf. \eqref{eq_def_Ahat_single_part}, and initial distribution $\mu_0$. The identity can be extended to all $h \in C([0,1])$ by standard arguments, see e.g. \cite[Section 6]{EK93}. 
In analogy to known duality properties of the microscopic particle system \cite{GKR07,CGR21} this can for example be used to get closed evolution equations for moments. Due to size-biasing $h(z)=z$ describes the expected second moment of the mass distribution and we get
\[
\frac{d}{dt} \E_{\mu_0} [\mu_t (z)] =\mathbf{E}_{\mu_0} [A_\theta \mathrm{Id}(
Z(t))]=\mathbf{E}_{\mu_0} \Big[ 2\big( 1-(1+\theta )Z_t\big)\Big] =2-2(1+\theta
)\E_{\mu_0} [\mu_t (Z)]\ .
\]
This has an exponential solution which converges to the stationary point $\frac{1}{1+\theta}$, the expected second moment of the GEM$(\theta )$ distribution. 
Dualities of this form were previously considered in \cite{DH82} and \cite{DK99}, see also \cite[Section 6]{EK93} for a summary. In Proposition \ref{prop_dLinf_duality} we will see that for $dL\to\infty$ dualities of the form \eqref{eq_dual2} extend directly to nonlinear test functions $H(\mu )$, due to the absence of an interaction part in the generator $\genInf$, cf. \eqref{eq_def_G}.
In the present case it may be possible to establish higher dimensional dual processes evaluated at $Bh$ rather than $h$, with  Fleming-Viot-type resampling. This is not relevant for the aim of this paper but could be an interesting question for future studies.

We stress once more the difference in point of view: whereas in previous works, see \cite{EK93} and references therein, the dual particles encode the position of clusters on the underlying lattice, in our size-biased approach the state of dual particles characterises the fragmentation of mass in a given configuration/partition. 
This allows for observing the dynamics of macroscopic cluster size distributions, while tracking only a finite number of dual particles.

The duality in \eqref{eq_dual2} is also interesting from a computational point of view, as it allows to continuously track the expected behaviour of the infinite dimensional process using only a finite-dimensional diffusion, without running into any absorbing states as is the case when observing a fixed set of lattice sites. 
A simple example is the second moment of cluster sizes in the Poisson-Dirichlet diffusion
at time $t$ which is given by $\mathbf{E}_{\mu_0}[Z(t)]$ as mentioned above.

\section{The diffusion limit in the case \texorpdfstring{$dL\to\infty$}{dL infinite}}\label{sec_dLinf}

The case of $dL\to \infty$ may be considered as an interpretation of the
Poisson-Dirichlet diffusion with infinite mutation rate $\theta = \infty$.
Clearly, this corresponds to an infinite drift towards zero in the single
particle operator $\mutFin$, cf. \eqref{eq_def_A_singlePart}. Thus, in order to see
non-trivial dynamics, we have to rescale time appropriately. 
Recalling \eqref{eq_IP_disc_approx_A}, we see that
\begin{align}
\frac{1}{dL}
    \mathfrak{L}_{L,N}H(\mu^{(\cdot)})(\eta)
    =
    \mu^{(\eta)}
    \Big(
     (h(0)-h(z)-z h'(z))
    \Big)+o(1)\,,
\end{align}
where $H(\mu)=\mu(h)$, $h\in C^3([0,1])$. The time-change also eradicates the interaction term in the corresponding limiting measure-valued process. We are left with a process that pushes mass (deterministically) from the interval $(0,1]$ onto zero.
Hence, mass will not accumulate on the macroscopic scale, instead we need to consider an appropriate mesoscopic scale to see the actual dynamics of the fast mixing mechanism.

In \cite{JCG19} it was proven that, at stationarity, mass accumulates on the mesoscopic scale of order $d^{-1}$, when $\rho\in (0,\infty)$, cf. \eqref{eq_exp_limi_JCG19}. Thus the embedding of particle configurations into $\mathcal{M}_1(\R_+)$ via \eqref{eq_def_map_mv_dLinf} with $\Hat{\mu}^{(\eta)}= \sum_{x=1}^L \frac{\eta_x}{N} \delta_{\frac{dL}{N} \eta_x}$
is an appropriate a-priori choice.\footnote{For the case $\rho\in (0,\infty)$,
also the choice $\delta_{d \eta_x}$ is appropriate and leads to a $\rho$
dependent limit. However, for the boundary cases $\rho \in \{0,\infty\}$ the
given choice turns out to be the correct one.}
In order to take particle configurations with mass lying on larger scales than $N/(dL)$ into account, we will consider probability measures $\mathcal{M}_{1}(\overline{\R}_+)$ on the one-point compactification, instead of restricting ourselves to the positive real line. 
We equip $\mathcal{M}_{1}(\overline{\R}_+)$ with the topology induced by weak convergence, thus, $\mathcal{M}_{1}(\overline{\R}_+)$ is compact. 

\subsection{Deriving the diffusion limit}

Once more we rely on the Trotter-Kurtz approximation to conclude the scaling limit in Theorem~\ref{theo_dL_inf}.
We follow the same steps as in Section~\ref{sec_dLfin}, carried out below for
completeness.

\begin{proposition}\label{prop_dLinf_multi_particle}
Let $H\in\mathcal{D}({\genInf})$, then 
\begin{align}
\lim_{N/L\to \rho} \sup_{\eta\in \Omega_{L,N}}
    \left|\frac{1}{dL}\mathfrak{L}_{L,N}H(\Hat{\mu}^{(\cdot)})(\eta) -
{\genInf}H (\Hat{\mu}^{(\eta)}) \right|=0\,.
\end{align}\end{proposition}

Below, we will show that the interaction term of the limiting measure-valued process indeed vanishes. First, we only consider test functions of the form $\mu\mapsto \mu(h)$.

\begin{lemma}\label{lem_dLinf_single_particle}
Let $\rho \in (0,\infty)$ and $H(\mu)=\mu(h)$, with $h\in \mathcal{D}({\mutInf})$ \eqref{eq_def_adomain}. Then
\begin{align}
\lim_{N/L\to \rho}
    \sup_{\eta\in \Omega_{L,N}}
    \left|\frac{1}{dL}\mathfrak{L}_{L,N}H(\Hat{\mu}^{(\cdot)})(\eta) -
\Hat{\mu}^{(\eta)}({\mutInf}h) \right|=0\,,
\end{align}
where ${\mutInf}$ is the single-particle generator defined in  \eqref{eq_def_Ahat_single_part}.
\end{lemma}

\begin{proof}
    Without loss of generality, we assume that $h\in C_c^3(\R_+)$.
    For simplicity of notation we will write $p_x=dL\tfrac{\eta_x}{N}$.
    Following the same steps as in the proof of  Lemma~\ref{lem_conv_single_part_dL_fin}, we have
    \begin{align}\label{eq_dLinf_taylor_single_part}
    \frac{1}{dL}
    \mathfrak{L}_{L,N}H(\Hat{\mu}^{(\cdot)})(\eta)
    &= 
    \frac{1}{2N^2}
    \sum_{\substack{x,y=1 \\ x\neq y}}^L 
    \eta_x \eta_y
    \big( 
    \widetilde{h}''(p_x)+\widetilde{h}''(p_y)
    \big)\\
    &\qquad+
    \frac{1}{L\,N}
    \sum_{\substack{x,y=1 \\ x\neq y}}^L 
    \eta_x
    \big(\widetilde{h}'(p_y)
    -\widetilde{h}'(p_x)\big)
    +o(1)\,, \nonumber 
\end{align}
where we gained an additional factor $(dL)^{-1}$ by rewriting $\tfrac{\eta_x}{N}h(d\eta_x)=\tfrac{1}{dL}\Tilde{h}(p_x)$.
Here we used again the fact that second-order terms in the second sum have a vanishing contribution and first-order terms in the first sum cancel exactly, cf. proof of Lemma~\ref{lem_conv_single_part_dL_fin}.
Display \eqref{eq_dLinf_taylor_single_part} can be written as
\begin{align}\label{eq_disc_dLinf_dyn_PD_diff_like_gen}
\frac{1}{dL}
    \mathfrak{L}_{L,N}H(\mu^{(\cdot)})(\eta)
    &=
    \sum_{\substack{x=1}}^{L}\frac{\eta_x}{N} \Big(
        1-\frac{p_x}{dL}
        \Big)
        \Tilde{h}''(p_x)\\
    &\qquad +
        \left[
        \frac{1}{L}\sum_{y=1}^{L} \Tilde{h}'(p_y)
        -
        \sum_{x=1}^{L} \frac{\eta_x}{N}
        \Tilde{h}'(p_x)
        \right]
        +o(1)\,. \nonumber
\end{align}
Using again Lemma~\ref{lem_riemann_partition_h0}, we have
$\frac{1}{L}\sum_{y=1}^{L} \Tilde{h}'(p_y)= \Tilde{h}'(0)+o(1)$. Hence, 
\begin{align}
\frac{1}{dL}
    \mathfrak{L}_{L,N}H(\mu^{(\cdot)})(\eta)
    =
    \Hat{\mu}^{(\eta)}\left( 
        \Tilde{h}''(p)
    +
        [
        \Tilde{h}'(0)
        -
        \Tilde{h}'(p)
        ]\right)
        +o(1)
        = \Hat{\mu}^{(\eta)}\big( 
        {\mutInf}h
        \big)+o(1)\,,
\end{align}
where we additionally used the fact that $\Tilde{h}''$ is bounded and of compact support, thus, $\tfrac{p}{dL}\Tilde{h}''(p)$ vanishes in the thermodynamic limit because $dL\to \infty$. 
\end{proof}

\begin{remark}\label{rem_dL_ext_generator}
Considering the state space $\mathcal{M}_{\leq 1}(\R_+)$ (positive measures on $\R_+$ with mass $\leq 1$) instead of $\mathcal{M}_1(\overline{\R}_+)$,
 \eqref{eq_disc_dLinf_dyn_PD_diff_like_gen} suggests that ${\mutInf}$ should act via
\begin{align}\label{eq_extended_gen_dL_inf}
    H \mapsto \big(\mu \mapsto \mu({\mutInf}h)+(1-\mu(1)) h(0)\big)\,,
\end{align}
on functions $H: \mathcal{M}_{\leq 1}(\R_+)\mapsto \R$ of the form $H(\mu)=\mu(h)$.
The extra term takes into account the transfer of mass from larger scales, i.e. above $N/(dL)$, which is pushed onto microscopic scales, cf. Corollary~\ref{cor_mass_process}. However, $\mu \mapsto \mu(1)$ is not a continuous function on $\mathcal{M}_{\leq 1}(\R_+)$. Instead the mass transport from larger scales is implicit in the generator ${\genInf}$, as we will see below.
\end{remark}

\begin{proof}[Proof of Proposition~\ref{prop_dLinf_multi_particle}]
    It suffices to consider functions $H\in \mathcal{D}({\genInf})$ of the form $H(\mu) = \mu(h_1)\cdots \mu(h_n)$, with $h_k\in C_c^3(\R)$.
    The interaction term of the operator ${\genInf}$ is again given by the second-order term of the following expansion
    \begin{align}\label{eq_dLinf_expansion_H_particle_shift}
    H(\Hat{\mu}_{\#} \eta^{x,y})
    &=
    H(\Hat{\mu}^{(\eta)}) \nonumber\\
    &\qquad+ \frac{1}{dL}\sum_{k=1}^n [\widetilde{h}_k(p_y+\tfrac{dL}{N})- \widetilde{h}_k(p_y)
    +
    \widetilde{h}_k(p_x-\tfrac{dL}{N})-\widetilde{h}_k(p_x)]\prod_{\substack{l=1\\l\neq k}}^n \Hat{\mu}^{(\eta)}(h_l)\nonumber \\
    &\qquad +\frac{1}{(dL)^2}
    \sum_{1\leq k <l \leq n}
    [\widetilde{h}_k(p_y+\tfrac{dL}{N})- \widetilde{h}_k(p_y)
    +
    \widetilde{h}_k(p_x-\tfrac{dL}{N})-\widetilde{h}_k(p_x)]\\
    &\qquad \qquad \times 
    [\widetilde{h}_l(p_y+\tfrac{dL}{N})- \widetilde{h}_l(p_y)
    +
    \widetilde{h}_l(p_x-\tfrac{dL}{N})-\widetilde{h}_l(p_x)]
    \prod_{\substack{j=1\\j\neq k,l}}^n \Hat{\mu}^{(\eta)}(h_j)\nonumber\\
    &\qquad+ \frac{1}{(dL)^3}r(\eta)\,. \nonumber
\end{align}
In the first order term each summand can be treated individually using Lemma~\ref{lem_dLinf_single_particle}, it only remains to check that both second-order term and remainder have no contribution.

Using a first-order Taylor expansion yields the following bound
\begin{align*}
    &\left|\widetilde{h}(p_y+\tfrac{dL}{N})- \widetilde{h}(p_y)
    +
    \widetilde{h}(p_x-\tfrac{dL}{N})-\widetilde{h}(p_x)
    \right|\\
    &\qquad =
    \tfrac{dL}{N}\left|\Tilde{h}'(p_y) - \Tilde{h}'(p_x)+\frac{1}{2}\tfrac{dL}{N}(\Tilde{h}''(\xi_y)+\Tilde{h}''(\xi_x))\right|\\
    & \qquad \leq 2\tfrac{dL}{N} (\|\Tilde{h}'\|_{\infty}+\|\Tilde{h}''\|_{\infty}),
\end{align*}
where $\xi_x,\xi_y\in [0,dN]$ corresponding to the associated remainder term.
Hence, the second-order term in \eqref{eq_dLinf_expansion_H_particle_shift}, after applying $\tfrac{1}{dL}\mathfrak{L}_{L,N}$, is upper bounded (up to a constant) by 
\begin{align*}
    \frac{1}{d L}
    \sum_{1\leq k <l \leq n}
    \frac{1}{N^2}
    \sum_{\substack{x,y=1\\x\neq y}}^L {\eta_x} ({d}+\eta_y)
    (\|h_k'\|_{\infty}+\|h_k''\|_{\infty})
    (\|h_l'\|_{\infty}+\|h_l''\|_{\infty})
    \prod_{\substack{j=1\\j\neq k,l}}^n \|h_j\|_\infty
    \,,
\end{align*}
which vanishes as $dL\to\infty$. For the same reason, higher order terms in the expansion \eqref{eq_dLinf_expansion_H_particle_shift} have no contribution either.
\end{proof}

\sloppy The closure of $({\mutInf},\mathcal{D}({\mutInf}))$ generates a Feller semigroup on
$\overline{\R}_+$, thus, the closure of
$({\genInf},\mathcal{D}({\genInf}))$ generates a Fleming-Viot process on
the compact space $\mathcal{M}_1(\overline{\R}_+)$, with trajectories in
$C([0,\infty),\mathcal{M}_1(\overline{\R}_+))$, in the absence of interaction \cite[Theorem 2.3]{EK87}.
We now have everything at hand to prove our second main result, Theorem~\ref{theo_dL_inf}.

\begin{proof}[Proof of Theorem~\ref{theo_dL_inf}]
Once more, we apply \cite[Theorem IV.2.11]{EK_book} together with Proposition~\ref{prop_dLinf_multi_particle}, which immediately concludes the desired convergence \eqref{eq_theo_dLinf_conv} in $D([0,\infty),\mathcal{M}_1(\Bar{\R}_+))$. Note that any fixed initial condition $\mu\in \mathcal{M}_1(\Bar{\R}_+)$ can be approximated by particle configurations using the embedding $\Hat{\mu}^{(\cdot)}$, cf. Lemma~\ref{lem_particle_config_approx_dLinf}.
This completes the proof.
\end{proof}

\subsection{Duality and the hydrodynamic limit}

The absence of interaction in ${\genInf}$ leads to a deterministic evolution
of $( \Hat{\mu}_t(h))_{t\geq 0}$, $h\in \mathcal{D}({\mutInf})$. 
This is a consequence of Dynkin's formula as derived in \eqref{eq_sdes}, 
since the process 
solves the ODE $d \Hat{\mu}_t(h) = \Hat{\mu}_t({\mutInf}h)\,dt$.
Hence, the evolution of $( \Hat{\mu}_t(h))_{t\geq 0}$ can be described by a
single particle evolving according to the process generated by $ {\mutInf}$,
averaged over its initial condition $ \Hat{\mu}_{0}$. 
See also the duality mentioned in 
\eqref{eq_dual} which we extend in the next result. 
Unlike in the case of $dL\to \theta<\infty$, we can fully characterise the semigroup of ${\genInf}$ by considering only the evolution w.r.t. the single particle generator ${\mutInf}$. 

\begin{proposition}\label{prop_dLinf_duality}
Let $g\in C(\overline{\R}_+^n)$ and define $G(\mu):=\mu^{\otimes n}(g)$, then for any $\Hat{\mu}_0\in \mathcal{M}_1(\overline{\R}_+)$
\begin{align}\label{c}
    G\left(\Hat{\mu}_t\right)
    =
    \mathbf{E}_{\Hat{\mu}_0^{\otimes n}}\left[g(\Hat{Z}(t))\right]\,,\quad \forall t\geq 0\,,
\end{align}
where $ ( \Hat{\mu}_{t})_{ t \geqslant 0}$ evolves w.r.t. $\genInf$
and
 $(\Hat{Z}(t))_{t\geq 0}$ is the process consisting of $n$ independent copies generated by the single-particle generator ${\mutInf}$, cf. \eqref{eq_def_Ahat_single_part}.
In particular, for $n=1$  we have $\Hat{\mu}_t = Law(\Hat{Z}(t))$ whenever the initial conditions agree in the sense that $\Hat{Z}(0)\sim \Hat{\mu}_0$.
\end{proposition}


\begin{proof}

Following precisely the same steps as in \cite[Section 6]{EK93} using the resolvent operator, one can 
conclude the duality
\begin{equation}\label{eq_duality_supp}
\begin{aligned}
\E_{ \hat{\mu}_{0}} \left[  G\left(\Hat{\mu}_t\right)\right]
    =
    \mathbf{E}_{\Hat{\mu}_0^{\otimes n}}\left[g(\Hat{Z}(t))\right]\,,\quad
\forall t\geq 0\,,
\end{aligned}
\end{equation}
with $g\in C(\overline{\R}_+^n)$ and $G(\mu):=\mu^{\otimes n}(g)$. 
This is essentially a direct consequence of the absence of an interaction term in the generator $\genInf$, cf. \eqref{eq_def_G}, implying for
$H(\mu)=\mu(h_1)\cdots
\mu(h_n)$, $h_i\in \mathcal{D}(\mutInf)$,
\[
\genInf H(\mu )=
\sum_{1\leq k \leq n }
    \mu({\mutInf}h_k)
    \prod_{\substack{l=1 \\ l\neq k}}^n \mu(h_l)\,.
\]

Now, let us consider the case $n=1$ for which the identity reads
\begin{equation}\label{eq_detevol}
\begin{aligned}
\Hat{\mu}_t(h)
=
\E_{ \hat{\mu}_{0}} \left[  \Hat{\mu}_t(h)\right]
    =
    \mathbf{E}_{\Hat{\mu}_0}\left[h(\Hat{Z}(t))\right]\,,\quad
\forall t\geq 0\,,\  h \in C ( \overline{\R}_{+}) \,,
\end{aligned}
\end{equation}
where we additionally used the fact that $ \hat{\mu}_{t}(h)$ is deterministic,
cf. \eqref{eq_sdes}.
In particular, \eqref{eq_detevol} implies that $\Hat{\mu}_t = Law(\Hat{Z}(t))$
and the measure-valued evolution $(\Hat{\mu}_t)_{ t \geqslant 0}$ is indeed
deterministic. 
Hence, the expected value on the left-hand side of \eqref{eq_duality_supp} has
no affect and can be dropped.
\end{proof}


The duality result in Proposition~\ref{prop_dLinf_duality}, and equivalently
Theorem~\ref{theo_dL_inf}, can be interpreted in the sense of a hydrodynamic limit.

\begin{proposition}[Hydrodynamic limit]\label{prop_fokker_plank}
Consider the process $(\Hat{\mu}_t)_{t\geq 0}$ generated by ${\genInf}$ with initial data $\Hat{\mu}_0\in \mathcal{M}_1(\overline{\R}_+)$. 
Then for every $t>0$, $\Hat{\mu}_t$ has a Lebesgue-density $f(t,\cdot)$ on $\R_+$. 
The evolution of the density $(f(t,\cdot))_{t> 0}$  solves 
\begin{align}\label{eq_fokker_plank}
\begin{cases}
    \partial_t f(t,z) = z\,\partial_z^2f(t,z)
    +z\, \partial_z f(t,z)\\
    \lim_{z\to 0}f(t,z)=1
\end{cases} 
\end{align}
with $\lim_{t\to 0}\int_0^\infty h(z) f(t,z) \, dz=\Hat{\mu}_0(h)
$ for every $ h\in C_c(\R_+)$. 
\end{proposition}

\begin{remark}
The diffusion part of the generator ${\mutInf}$, given by \begin{align}\label{eq_cir_model}
    {\mutInf}_Dh(z):=z h''(z)+ (2-z)h'(z)\,,
\end{align}
generates the so called Cox-Ingersoll-Ross model \cite{DeMarco11}, which is a well studied
diffusion process in mathematical finance and population genetics.
\end{remark}

\begin{proof}
Consider first the case for an initial condition that has no atom at infinity,
i.e. $\|f_0\|_{L^1(\R_+)}=1$, and 
let $z_0\in [0,\infty)$. 
The Cox-Ingersoll-Ross model generated by ${\mutInf}_D$ is known to have a
density $g(t,\cdot |z_0)$ for any positive time and initial data \cite[Display
below (3.2)]{DeMarco11}. In fact, it is explicitly given and for $z_0=0$ it evaluates to 
\begin{align}
    g(t,z|0)
    =
    \frac{z}{(2\ell_t)^2} e^{-z(2\ell_t)^{-1}}\quad\mbox{with}\quad\ell_t:= \tfrac{1}{2}(1-e^{-t})\,.
\end{align}
Furthermore, for any $t>0$ we have $g_t(\cdot| z_0)\big|_{(0,\infty)}\in C^\infty((0,\infty))$ \cite[Proposition~3.2]{DeMarco11}.
The resetting mechanism is given by a Poisson jump process, thus, \cite[Theorem~1]{APZ13} guarantees that also the process $(\Hat{Z}(t))_{t\geq 0}$ generated by ${\mutInf}$ has a density which is given by 
\begin{align}
    f(t,z|z_0)=
    e^{-t}g(t,z|z_0)
    +
    \int_0^t e^{-s} g(s,z|0) \,ds\,.
\end{align}
We note that $f(t,\cdot|z_0)$ inherits the regularity properties of
$g(t,\cdot|z_0)$ on $(0,\infty)$. This follows from the change of variable $r =
z/(2\ell_s)$ with $\tfrac{dr}{ds}=-\tfrac{z}{(2\ell_s^2)}e^{-s}$ which yields for every $z>0$
\begin{align}\label{eq_density_jump_diff} 
    f(t,z|z_0)=
    e^{-t}g(t,z|z_0)
    +
    \int_{\frac{z}{2\ell_t}}^\infty
    e^{-r}
    \, dr
    =
     e^{-t}g(t,z|z_0)
     +e^{-z ( 2 \ell_t )^{-1}}
    \,.
\end{align}
Thus, $ f(t,\cdot |z_0) \in C^\infty((0,\infty))$.

Now, it is only left to verify that $(f(t,\cdot))_{t\geq 0}$ indeed solves the given PDE. Using integration by parts, we see that for any $h\in C_c^2(\R_+)$, we have 
\begin{align}
    \mu_t({\mutInf}h) = 
    \int_0^\infty f(t,z)\,{\mutInf}h(z)\, dz
    =
    \int_0^\infty {\mutInf}^*f(t,z)\,h(z)\, dz
\end{align}
with the adjoint action defined as 
\begin{align}\label{eq_Ahat_adjoint}
        {\mutInf}^*f(z):=
        zf''(z) +zf'(z)
        + \delta_0(z) \big(1-f(z)\big)\,.
\end{align}
Hence, the density $f=\big(f(t,\cdot)\big|_{(0,\infty)}\big)_{t\geq 0}$ in \eqref{eq_density_jump_diff} solves the PDE 
\begin{align}
    \partial_t f ={\mutInf}^*f\,, \quad f(0,\cdot)=\delta_{z_0}\,.
\end{align}
It is easy to see from \eqref{eq_density_jump_diff} that  $\lim_{z\to
0}f(t,z)=1$ for any $t>0$ since  $g(t,0|z_0)=0$ for all $z_{0} \geq 0$, thus, the boundary term in \eqref{eq_Ahat_adjoint} vanishes and we are left with the PDE in the statement.

Consider now the case of $f_0=\delta_\infty$.
The only way $(\Hat{Z}(t))_{t\geq 0}$ can escape infinity is via the resetting mechanism.
Thus, \eqref{eq_density_jump_diff} still applies with $g(t,z|z_0)$ replaced by $\delta_\infty$, 
since $\mathbb{P}(\Hat{Z}(t)=\infty)$ is equal to the probability that the process has not jumped yet.
One can check that also $f(\cdot|\infty)$ solves the given PDE on $(0,\infty)$ with the correct boundary condition.
The result for arbitrary initial conditions now follows by integrating over the densities w.r.t. $\mu_0$ and Leibniz rule.
\end{proof}

Given the explicit form of the density \eqref{eq_density_jump_diff}, we can read of the evolution of the mass process $(\mu_t[\R_+])_{t\geq 0}$.

\begin{corollary}\label{cor_mass_process}
We have 
\begin{align}\label{eq_evolut_mass_proc}
    \Hat{\mu}_t[\R_+]=1-(1-\Hat{\mu}_0[\R_+])e^{-t}\,.
\end{align}
\end{corollary}

\begin{proof}
We integrate the PDE from Proposition~\ref{prop_fokker_plank} in space, which yields 
\begin{align}
    \partial_t \|f(t,\cdot)\|_{L^{1}(\R_+)}
    =
    \int_0^\infty \big( z\,\partial_z^2f(t,z)
    +z\, \partial_z f(t,z)\big)\, dz\,.
\end{align}
The right-hand side simplifies to the differential equation
\begin{align}\label{eq_ODE_mass_process}
     d\alpha_t = (1-\alpha_t)\, dt\,, \quad \alpha_0=\mu_0[\R_+]\,,
\end{align}
using integration by parts, its solution is given by \eqref{eq_evolut_mass_proc}.
\end{proof}

\begin{figure}[t]
\centering
\begin{subfigure}[b]{.32\linewidth}
\includegraphics[width=\linewidth]{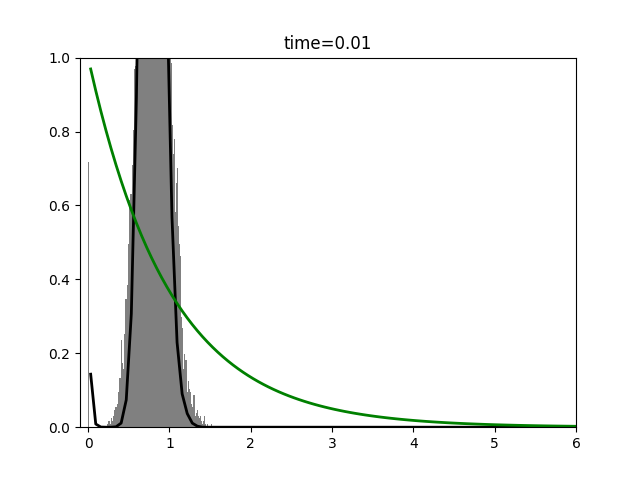}
\setcounter{subfigure}{2}%
\end{subfigure}
\begin{subfigure}[b]{.32\linewidth}
\includegraphics[width=\linewidth]{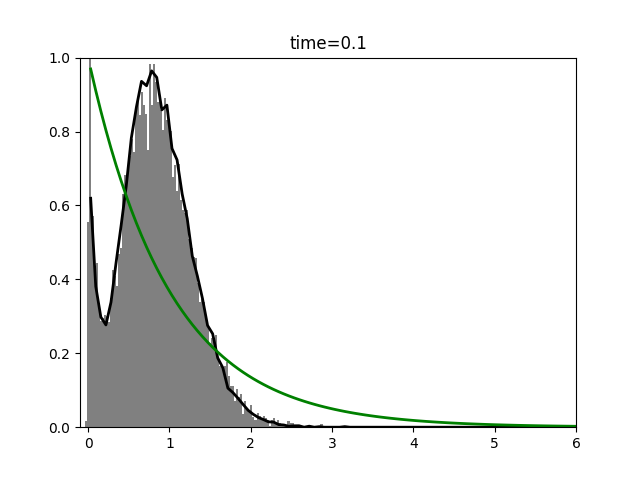}
\setcounter{subfigure}{2}%
\end{subfigure}
\begin{subfigure}[b]{.32\linewidth}
\includegraphics[width=\linewidth]{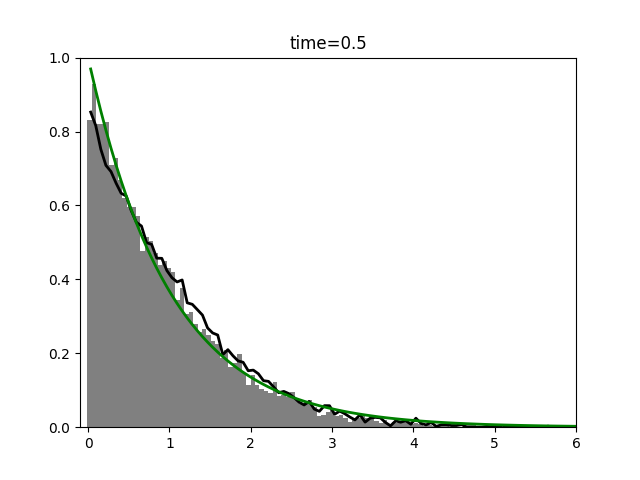}
\setcounter{subfigure}{2}%
\end{subfigure}
\caption{
Simulations for both the inclusion process ($N=L=1024$, $d= L^{-1/2}= \tfrac{1}{32}$) and the jump diffusion generated by ${\mutInf}$ agree in accordance with Proposition~\ref{prop_dLinf_duality}.
The black graph shows the density profile of embedded inclusion process \eqref{eq_def_map_mv_dLinf} (1000
samples), whereas
the grey histogram represents the density of the jump diffusion (10000 samples). 
Both profiles converge rapidly to the unit exponential density (green line), cf. Lemma~\ref{lem_prop_Ahat}.
We considered an initial condition $ \Hat{\mu}_0=  \delta_{ z}$,
$z \simeq \tfrac{25}{32}$. 
}
\end{figure}

Furthermore, we summarise invariance and exponential ergodicity of $(\Hat{\mu}_t)_{t\geq 0 }$ in the following lemma:
\begin{lemma}\label{lem_prop_Ahat}
The process $(\Hat{Z}(t))_{t\geq 0 }$ satisfies the following properties:
\begin{chosenEnum}
        \item The exponential distribution $\mathrm{Exp}(1)$ is the unique invariant probability measure.
        
        \item We have
        \begin{align}
            \left\| \text{Law}(\Hat{Z}(t))- \mathrm{Exp}(1)\right\|_{TV}\leq e^{- t}\,, \quad \forall t\geq 0\,.
        \end{align}
\end{chosenEnum}
\end{lemma}

\begin{proof}
The fact that the exponential distribution is invariant can be explicitly proven using integration by parts, but also follows directly from \cite[Corollary 1]{APZ13}. 
The exponential ergodicity is a consequence of \cite[Theorem 2]{APZ13}. The same result yields that the process is Harris recurrent, thus, the Exponential distribution is the unique invariant distribution.
\end{proof}

Lastly, it is an easy consequence that the point mass on the Exponential distribution is an absorbing point for the measure valued process.
\begin{corollary}\label{cor_exp_rever}
$\Hat{\mathbf{P}}= \delta_{\mathrm{Exp}(1)}$ is reversible w.r.t. the
dynamics induced by ${\genInf}$.
\end{corollary}

\begin{proof}
Let $H\in \mathcal{D}({\genInf})$ be of the form $H(\mu)=\mu(h_1) \cdots \mu(h_n)$. 
For simplicity we write $\nu := \mathrm{Exp}(1)\in \mathcal{M}_1(\R_+)$.
Then 
\begin{align}
    {\genInf}H(\nu) = \sum_{j=1}^n \nu({\mutInf}h_i) \prod_{l\neq j} \nu(h_l) =0\,,
\end{align}
as $\nu$ is invariant w.r.t. ${\mutInf}$, cf. Lemma~\ref{lem_prop_Ahat}. Thus, $\Hat{\mathbf{P}}(F {\genInf}H)=0= \Hat{\mathbf{P}}(H {\genInf}F)$.
\end{proof}

\subsection{A natural extension of the Poisson-Dirichlet diffusion}

In this section we prove Theorem~\ref{theo_theta_to_inf} which states that $\genFin$ (which is equivalent to the Poisson-Dirichlet diffusion) has a limit as $\theta\to \infty$ under appropriate rescaling.

\begin{proof}[Proof of Theorem~\ref{theo_theta_to_inf}]
The statement is once more a conclusion of the Trotter-Kurtz approximation. We state the essential steps for completeness.
Let $H\in \mathcal{D}({\genInf})$ be of the form $\mu \mapsto \mu(h_1) \cdots \mu(h_n)$, $h_k\in \mathcal{D}({\mutInf})$, and define $H_\theta \in \mathcal{D}(\genFin)$ by
\begin{align}
    E \ni \mu \mapsto \mu(h_1(\theta \, \cdot ))
    \cdots 
    \mu(h_n(\theta\,  \cdot ))\,,
\end{align}
where we interpreted the $h_k(\theta\, \cdot )$'s to be elements of $C^3([0,1])$. 
We have 
\begin{align}
    {\mutInf}h_k(\theta z ) -\frac{1}{\theta } {\mutFin} h_k(\theta \, \cdot)(z)
    = 
    z \big( \theta z\, h_k''(\theta z) +2\, h'_k(\theta z)\big)\,.
\end{align}
If $h_k$ is a constant function, the r.h.s. vanishes. On the other hand, if $h_k\in C^3_c(\R_+)$, we can write
\begin{align*}
    \sup_{\mu \in E}\big|
    \mu\big(({\mutInf}h_k)(\theta \cdot )\big)
    -
    \tfrac{1}{\theta }\mu\big({\mutFin} h_k(\theta \, \cdot)\big)
    \big|
    \leq &
    C_{h_k}\, \sup_{\mu \in E} \mu\big(Z\,  \mathds{1}_{\theta Z\in
\text{supp}(h_k)}\big)\\
    =& C_{h_k}\,\sup_{z\in [0,1]}
    \{z\, \mathds{1}_{\theta z \in \text{supp}(h_k)}\}\,,
\end{align*}
where $Z\sim \mu$ and $C_{h_k}$ is a finite constant, depending on $h_k$.
The right hand side vanishes as $\theta \to \infty$. 
Along the same lines, we can show that the interaction term of $\tfrac{1}{\theta}\genFin$ disappears. Overall we conclude
\begin{align}
\lim_{\theta \to \infty}
    \sup_{\mu \in E}\big| \tfrac{1}{\theta}\genFin H_\theta (\mu) - ({\genInf}H)(S_\theta \mu) \big|=0\,.
\end{align}
Again, the convergence of generators suffices to conclude weak convergence on the process level.
\end{proof}

\section{Discussion and outlook}\label{sec_overview}

We conclude this paper with a discussion of  boundary cases in the setting considered and outline future directions as well as work in progress.

Throughout the paper we assumed that $\rho=\lim_{N,L\to \infty} N/L\in
(0,\infty)$.
However, the derived scaling limits do not depend on the actual value of
\(\rho\), as we study the distribution of mass after renormalising by $N$.
As long as $N,L \to \infty$, our results extend to the boundary cases $\rho \in
\{0,\infty\}$, up to a
regime around $\rho =0$ in the case $dL \to \infty$. In this regime we see an interesting
transition of the clustering behaviour, cf. Lemma~\ref{lem_geometric_dist} below.

First consider $\theta<\infty$ and $d \to 0$, for which both cases $\rho \in \{0,\infty\}$ are covered by our proof. For $\rho=0$, i.e. $N\ll L$, this is clear intuitively, as an increasing number of empty sites does not affect the dynamics since the total diffusivity per particle is $dL\to\theta$.
On the other hand, if $\rho=\infty$ it may seem surprising that the number of sites $L$ does not play a role (as long they are divergent).
Here, the core lies in  Lemma~\ref{lem_particle_config_approx_simplex}, which states that for any thermodynamic limit, we can approximate configurations in $\overline{\nabla}$ (equivalently measures in $E$) by a sequence of particle configurations. Indeed, having a closer look at the proof of Lemma~\ref{lem_particle_config_approx_simplex}, we see it is only necessary that a macroscopic excess mass of order $\sim N$ can be distributed uniformly over sites such that it is not visible under the macroscopic rescaling $\tfrac{1}{N}$. This is always the case as we can put $\sim N/L$ particles on sites (which might itself diverge), however, under macroscopic rescaling we have $\sim \tfrac{1}{N}\tfrac{N}{L} \to 0$.

Now assume $\theta=\infty$. In the case $\rho\in (0,\infty)$, the results in the present paper, and also \cite{JCG19}, yield that a size-biased chosen chunk (at equilibrium) is approximately exponentially distributed with mean $\simeq \tfrac{N}{dL}$.
In fact, looking at the proof of Theorem~\ref{theo_dL_inf}, the result remains true as long as $dL/N \to 0$. This trivially holds when $\rho = \infty$, in which case cluster sizes live on the scale of order
\begin{align}
 \frac{N}{dL}    \gg \frac{1}{d}\,.
\end{align}
On the other hand, if $N/L\to 0$ we have no control over $dN$ (in contrast to $N/L \to \infty$, which implies $dN\to \infty$ since $dL\to \infty$).
Therefore, for $\rho=0$ the convergence in Theorem~\ref{theo_dL_inf} remains true only if $dL/N \to 0$, i.e. $d\ll N/L$. 
Assume on the other hand that $dL/N \to \gamma \in (0,\infty]$, then we don't expect any clustering of particles on diverging scales.
This is indeed the case, in fact, we see a finer structure emerging in the limit on scales of order one. 
Note that the following is independent of the underlying graph structure and holds more generally for irreducible and spatially homogeneous dynamics, cf. \cite{JCG19}.

\begin{lemma}\label{lem_geometric_dist}
    Assume that $N/L\to 0$, as $N,L\to\infty$, $d \to 0$ and $dL\to \infty$ such that $N/(dL) \to \gamma \in [0,\infty)$. Then 
    \begin{align}
        \lim_{N/L\to 0}\pi_{L,N}[\Tilde{\eta}_1\in \cdot ]
        =
        \mathrm{Geom}\left(\frac{1}{1+\gamma}\right)
                \,.
    \end{align}
    Here $\pi_{L,N}$ denotes the unique invariant distribution w.r.t. $\mathfrak{L}_{L,N}$.
\end{lemma}

Hence, for $\rho=0$ and $dL\to \infty$, there is a critical scaling $N\sim dL$ below which the equilibrium measure does not exhibit clustering of particles on diverging scales. 

\begin{figure}[t]
\centering
\begin{tikzpicture}
  \message{Phase transition of ice -> water -> steam^^J}
  \def\ymin{-0.5}
  \def\ymax{1.5}
  \def\xmin{0}
  \def\xmax{12}
  
 \node[below] at (0.5,0) {$\ll d$};
\node[below] at (2,0) {$\gamma\, d$};
\node[below] at (3.5,0) {$\gg d$};

 \node[below] at (0.5,1) {$\delta_{1,n}$};
 
  \node[below] at (2,1) {$\mathrm{Geom}(\tfrac{1}{1+\gamma})$};
  
   \node[below] at (7.5,1) {$\mathrm{Exp}(1)$ on scales $\frac{N}{dL}$};

  \draw[-,thick] (0,\ymin) -- (0,\ymax) node[below left=0] {};
    \draw[-,thick] (1,\ymin) -- (1,\ymax) node[below left=0] {};

    \draw[-,ultra thick] (3,\ymin-0.2) -- (3,\ymax+0.2) node[below left=0] {};
    \draw[-,dashed,thick] (4,\ymin) -- (4,0.5) node[below left=0] {};
  \draw[-,dashed,thick] (\xmax*0.9,\ymin) -- (\xmax*0.9,0.5) node[below left=0] {};
  \draw[-,thick] (\xmax,\ymin) -- (\xmax,\ymax) node[below left=0] {};

  \draw[-,thick] (\xmin,0) -- (\xmax,0) node[below left=0] {}; 
    \draw[-,thick] (\xmin,\ymax) -- (\xmax,\ymax) node[below left=0] {}; 
\draw [decorate,decoration={brace,amplitude=5pt,mirror,raise=4ex}]
  (4+0.1,0) -- (\xmax*0.9-0.1,0) node[midway,yshift=-3em]{$\frac{N}{L}\to \rho \in (0,\infty)$};
 \draw [decorate,decoration={brace,amplitude=5pt,mirror,raise=4ex}]
  (0,0) -- (4-0.1,0) node[midway,yshift=-3em]{$\rho=0$};
\draw [decorate,decoration={brace,amplitude=5pt,mirror,raise=4ex}]
  (\xmax*0.9+0.1,0) -- (\xmax,0) node[midway,yshift=-3em]{$\rho =\infty$};
\end{tikzpicture}
\caption{
Graphical summary of the clustering of particles at equilibrium for the inclusion process when $dL\to\infty$.
The distributions displayed describe the first size-biased marginal $\Tilde{\eta}_1$ on the appropriate scale.
Note particularly the transition from diverging scales to scales of order $1$,
when moving from the regime $\rho \gg d$ into $\rho \sim \gamma d$.
}
\end{figure}
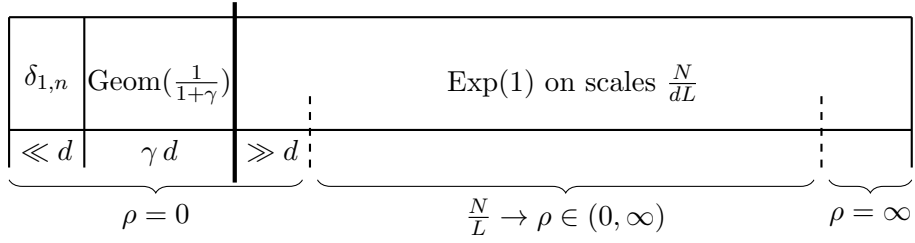

Before proving the above lemma, we require some notation and representations, which can be found in more detail in \cite{JCG19}.  
Recall that $\pi_{L,N}$ denotes the unique invariant distribution w.r.t. $\mathfrak{L}_{L,N}$ supported on $\Omega_{L,N}$, which is  given explicitly  by 
\begin{align}
    \pi_{L,N}[d\eta]=\frac{1}{Z_{L,N}}\prod_{x=1}^L
    w_L(\eta_x)\, d\eta\,,
\end{align}
where $d\eta$ denotes the counting measure, $Z_{L,N}$ the appropriate partition function and $w_L$ describing weights of the form 
\begin{align}\label{eq_weights}
    w_L(n)=
    \frac{\Gamma(n+d)}{n! \Gamma(d)}\,,
\end{align}
arising from the choice of transition rates in $\mathfrak{L}_{L,N}$.
Moreover, the partition function can be explicitly written in terms of 
\begin{align}\label{eq_pf}
    Z_{L,N}=\frac{\Gamma(N+dL)}{N!\Gamma(dL)}\,.
\end{align}

\begin{proof}[Proof of Lemma~\ref{lem_geometric_dist}]
Recall from \cite[Display (14)]{JCG19} that 
\begin{align}
    \pi_{L,N}[\Tilde{\eta}_1 =n]
    =
    \frac{L}{N} n\, w_{L}(n) \frac{Z_{L-1,N-n}}{Z_{L,N}}\,,
\end{align}
which is equal to zero if $n=0$. Thus, without loss of generality let $n>0$.
We replace the terms in the previous display with the corresponding expressions
in \eqref{eq_weights} and \eqref{eq_pf}, which yields 
\begin{align}
    \pi_{L,N}[\Tilde{\eta}_1 =n]
    &\simeq 
    \frac{dL}{N} 
    \frac{\Gamma(n+d)}{(n-1)!\, d \Gamma(d)}
    \frac{\Gamma(N-n+dL)}{\Gamma(N+dL)} N^n
    \,.
\end{align}
We analyse the remaining terms individually and conclude
\begin{align}
    \pi_{L,N}[\Tilde{\eta}_1 =n]
    &\simeq 
    \frac{dL}{N} 
    \bigg(\frac{\frac{N}{dL}}{1+\frac{N}{dL}}\bigg)^n
    \to 
    \frac{\gamma^{n-1}}{(1+\gamma)^n}\,, \quad \text{if }\gamma\in [0,\infty)\,.
\end{align}
This finishes the proof.
\end{proof}

Our approach can further be adapted, with minimal changes in Lemma~\ref{lem_riemann_partition_h0}, to cover the situation of fixed $L$ and $d =d_N \to 0$ as $N\to\infty$. In this limit, the full mass fraction condenses on a single cluster site in a typical stationary configuration. This scaling has been considered in \cite{BDG17,KS21} with additional assusmption $d_N \log N\to 0$ to study the metastable dynamics of the cluster location on a large time scale, and in \cite{GRV13} under the additional assumption $d_N N \to\infty$ to study emergence of the single cluster site on a slower time scale. The interaction and metastable motion of several spatially separated clusters on arbitrary finite lattices has been understood in full detail by now (see \cite{kim2023hierarchical} and references therein). In our model clusters can exchange particles directly due to the complete graph geometry. Our results in both regimes imply that their rescaled sizes fluctuate diffusively due to mass exchanges on the time scale $1/d_L$. (Note that the generator \eqref{eq_def_ip_dynamics} has already been speeded up by a factor $L-1$ for complete graph dynamics). Only in the case $d_L L\to 0$ does a single stable cluster emerge, whose location should exhibit a uniform jump process on a time scale $1/d_L$. This is not covered by our results and could be approached by a combination of techniques in \cite{BDG17,KS21} and \cite{armendariz2017metastability}.\\

In summary, our approach not only fully determines the clustering of particles on diverging scales in the inclusion process, but also describes the corresponding dynamics of the limiting Markov process in case of the complete graph.
In this paper we focused on a (joint) thermodynamic limit $N,L\to\infty$ with
$N/L\to\rho$ which is a natural approach in the context of interacting particle
systems. It is equally possible to consider a two-step limit
$N\to\infty$ towards a Wright--Fisher diffusion on a finite dimensional simplex,
then taking $L\to\infty$ in a second step, arriving at the same scaling limits.
This approach is often used in the context of population genetics
\cite{EK81,CBERS17,RW09}. 
Following the discussion above, this should be possible
with the same size-biased techniques used in this paper, using  a suitable scaling of $d$.
\\

A natural next step, in view of the hydrodynamic limit  Proposition~\ref{prop_fokker_plank}, when $dL\to\infty$, is to study fluctuations around the equilibrium.
A second moment calculation w.r.t. the stationary measure yields 
\begin{align}
        \pi_{L,N}\left(\big(\Hat{\mu}^{(\eta)}(h) - \mathrm{Exp}(1)(h) \big)^2\right)\simeq \frac{1}{dL} \pi_{L,N}\big( \Tilde{\eta}_1\tfrac{dL}{N} h(\Tilde{\eta}_1\tfrac{dL}{N})^2 \big) \to 0\,,
\end{align}
for any $h\in C_b(\R_+)$. Hence, in order to see a non-trivial limit, we should investigate fluctuations of order $\sqrt{dL}$ by studying the limiting behaviour of
\begin{align}
    \sqrt{dL}\left((\Hat{\mu}_{\#} \eta^{(L,N)}(\tfrac{t}{dL}))(h) - \mathrm{Exp}(1)(h) \right)\,, \quad t\geq 0\,,
\end{align}
where $(\eta^{(L,N)}(t))_{t\geq 0}$ denotes the inclusion process generated by $\mathfrak{L}_{L,N}$. Note that we slowed down time of the process, as indicated in Theorem~\ref{theo_dL_inf}.
Due to decoupling of the size-biased marginals, the fluctuations are expected to be Gaussian.\\

    
    

Our approach should be robust towards perturbation of transition rates, as we do not require the explicit form of the partition function of the canonical distribution, cf. \eqref{eq_pf}. However, the compact form of the limit dynamics and in particular the duality \eqref{eq_dual} are not expected to extend to more general models. 
Throughout the paper, we have focused on the one parameter family of Poisson-Dirichlet diffusions. There exists a two-parameter extension of the process, which was introduced in \cite{Pe09}. This process has gained a lot of attention over the past years \cite{RW09,FSWX11,Et14,CBERS17}, just to name a few.
It would be interesting to investigate the size-biased approach in this setting. 
To the best of the authors' knowledge, the two parameter process has only been
studied when fixing finitely many locations/sites and observing the evolution of mass on them, see for example \cite{FSRW21} and also the discussion in Section~\ref{sec_cf_sb_fix_site}.

Furthermore, it would be interesting to investigate diffusion limits of the generalised version of the inclusion process with non-trivial bulk, studied in \cite{CGG21}.
Numerical simulations and heuristic arguments suggest that the macroscopic phase evolves under the dynamics described in Theorem~\ref{theo_dL_fin}.
At the same time one can observe a transfer of mass between the bulk and the condensate whose evolution is described by a system of ODEs, similar to Corollary~\ref{cor_mass_process}.




%% file: appendix.tex
\section{Embeddings and approximations of particle configurations}

First, we show that size-biased probability measures $E$ are isomorphic to the Kingman simplex.

\begin{lemma}\label{lem_isomorphism_simplex_measures}
The map $\mu^{(\cdot)}: \overline{\nabla} \to E$, cf. \eqref{eq_embedding_kingman_to_measure}, is an isomorphism. 
\end{lemma}
\begin{proof}
First note that surjectivity is trivial due to the definition of $E$.
Now, consider $p,q\in \overline{\nabla}$ such that $p\neq q$. Then there exists
an index $i\in\N$ such that $p_i\neq q_i$ and $p_j=q_j$ for all $j<i$, without loss of generality assume $p_i>q_i$. Then 
\begin{align}
    \mu^{(p)}([p_i,1])\geq \sum_{j=1}^i p_j > \sum_{j=1}^{i-1} q_j =\mu^{(q)}([p_i,1])\,.
\end{align}
thus, $\mu^{(p)}\neq \mu^{(q)}$. 

In order to show that the map $\mu^{(\cdot)}$ is continuous, consider a sequence of partitions $(p^{(n)})_{n\in \N}$ converging to $p$ in $\overline{\nabla}$. Then for every $h\in C([0,1])$
(uniformly in $n$)
    \begin{align}\label{eq_cont_of_testing}
        \Big|
        \sum_{i=M}^\infty p_i^{(n)} (h(p_i^{(n)})-h(0))
        \Big|
        \leq 
        \sup_{0\leq z\leq \frac{1}{M}} |h(z)-h(0)|\to 0\,, \quad \text{as }M\to \infty\,,
    \end{align}
    where we used the fact that $p_i\leq \tfrac{1}{i}$ for any $i\in \N$. This implies in particular $\mu^{(p^{(n)})}\weakconv \mu^{(p)}$, recall that
    $\mu^{(p)}(h)=h(0)+ \sum_{i=1}^\infty p_i (h(p_i)-h(0))$.
    
Continuity of the inverse is now immediate: let $(\mu_n)_{n\in \N}$ be a sequence in $E$ weakly converging to $\mu\in \mathcal{M}_1([0,1])$. Then we can identify each $\mu_n$ with a unique $p^{(n)}$ satisfying $\mu^{(p^{(n)})}=\mu_n$. Due to compactness of $\overline{\nabla}$, it suffices to consider convergent subsequences, say $(p^{(n_j)})_{j\in \N}$ with limit $p$. Thus, by assumption and continuity of $\mu^{(\cdot)}$
\begin{align}
    \mu_{\#} p^{(n_j)} \weakconv \mu^{(p)} = \mu\,,
\end{align}
which particularly implies that $\mu\in E$. This implies that each accumulation point must agree with $(\mu^{(\cdot)})^{-1}(\mu)=p$.
\end{proof}

\begin{lemma}\label{lem_riemann_partition_h0}
Let $h\in C(\R_+)$ and $\rho\in[0,\infty)$. Then for any $\zeta_L \to 0$ 
\begin{align}
    \lim_{N/L \to \rho}
    \sup_{\eta \in \Omega_{L,N}}
    \Big|
    \frac{1}{L}\sum_{x=1}^L  h(\zeta_L\eta_x) -h(0)
    \Big|= 0\,.
\end{align}
\end{lemma}

\begin{proof}
Let $\varepsilon>0$ and $\eta\in \Omega_{L,N}$, then
\begin{align*}
\Big|
    \frac{1}{L}\sum_{x=1}^L  h(\zeta_L\eta_x)
    -h(0)\Big|
    &\leq \Big|
    \frac{1}{L}\sum_{x=1}^L \big( h(\zeta_L \eta_x) -h(0)\big)
    \mathds{1}_{\zeta_L\eta_x > \varepsilon }\Big|\\
    &\qquad +\Big|
    \frac{1}{L}\sum_{x=1}^L
    \mathds{1}_{\zeta_L\eta_x \leq \varepsilon }
    \big( h(\zeta_L \eta_x) -h(0)\big)\Big|\\
    &\leq
    2\|h\|_{\infty}
    \frac{1}{L}\sum_{x=1}^L
    \mathds{1}_{\zeta_L\eta_x > \varepsilon }
    +
    \frac{1}{L}\sum_{x=1}^L  
    \mathds{1}_{\zeta_L\eta_x \leq \varepsilon }
    \big| h(\zeta_L \eta_x) -h(0)\big|\\
    &\leq 
    2\|h\|_{\infty}
    \frac{1}{L}\sum_{x=1}^L
    \mathds{1}_{\zeta_L\eta_x > \varepsilon }
    +
    \sup_{0\leq v\leq \varepsilon}\big| h(v) -h(0)\big|
    \frac{1}{L}\sum_{x=1}^L 
    \mathds{1}_{\zeta_L\eta_x \leq \varepsilon }\,.
\end{align*}
The first term on the r.h.s. vanishes because the number of sites satisfying
$\zeta_L \eta_x > \varepsilon$ is upper bounded by $\zeta_L\,N\,
\varepsilon^{-1}$ (otherwise the total mass exceeds $N$). Thus,
\begin{align}\label{eq_riemann_approx_crux}
    \frac{1}{L}\sum_{x=1}^L
    \mathds{1}_{\eta_x > \varepsilon \zeta_L^{-1}}
    \leq \frac{\zeta_L\, N }{\varepsilon L}\,.
\end{align}
The second term, on the other hand, is upper bounded by $\sup_{0\leq v\leq \varepsilon}\big| h(v) -h(0)\big|$, which vanishes in the small $\varepsilon$-limit.
Note that both upper bounds are uniform in $\Omega_{L,N}$. 
Now, taking first the thermodynamic limit $N/L \to \rho$ before taking $\varepsilon \to 0$, finishes the proof.
\end{proof}

\begin{remark}\label{rem_rho_inf_riemann_approx}
Note that Lemma~\ref{lem_riemann_partition_h0} remains true if $\rho=\infty$
with a choice $\zeta_L$ satisfying $\zeta_L L/N\to 0$, cf.~\eqref{eq_riemann_approx_crux}. 
In particular, the case $\zeta_L = dL/N$ is covered.
\end{remark}

Indeed, we can show that, independently of the thermodynamic limit taken, any element in $\overline{\nabla}$ can be approximated by particle configurations:
\begin{lemma}\label{lem_particle_config_approx_simplex}
    Let $N/L\to\rho\in [0,\infty]$, then
    for any $p\in \overline{\nabla}$ there exist $\eta^{(L,N)}\in \Omega_{L,N}$ such that 
    \begin{align}
        \frac{1}{N}\Hat{\eta}^{(L,N)} =\frac{1}{N}\eta^{(L,N)}\to p\,, \quad \text{in }\overline{\nabla}\,.
    \end{align}
\end{lemma}
\begin{proof}
    Consider $p\in \overline{\nabla}$ with $\|p\|_1=1-\gamma$, we then define $\Bar{\eta}^{(L,N)}_i:= \lfloor p_i\, N \rfloor$, for $i\in \{1,\ldots, L\}$. Hence, there are 
    \begin{align}
       M_{L,N}(p):=N- \sum_{x=1}^L \Bar{\eta}^{(L,N)}_x = \gamma\, N+ \sum_{i=1}^L(p_i N - \lfloor p_i\, N \rfloor) \leq \gamma\, N +L
    \end{align}
    particles to spare. Thus, defining $\eta^{(L,N)}\in \Omega_{L,N}$ via
    \begin{align}
        \eta^{(L,N)}_x
        :=
        \Bar{\eta}^{(L,N)}_x
        +\Big\lfloor \frac{M_{L,N}(p)}{L} \Big\rfloor
        +
        \mathds{1}_{x\leq ( M_{L,N}(p) \mod L ) }\,, 
        \quad 
        x\in \{ 1, \ldots , L\} \, ,
    \end{align}
    yields the desired approximation, since $M_{L,N}(p)/(N\,L) \to 0$, as $N,L\to\infty$.
\end{proof}

Similarly, the embedding via $\Hat{\mu}$, cf. \eqref{eq_def_map_mv_dLinf}, allows to approximate any probability measure on $\overline{\R}_+$ by particle configurations.

\begin{lemma}\label{lem_particle_config_approx_dLinf}
Let $ \rho \in [0, \infty]$ and  $d=d(L)$ such that 
\begin{equation*}
\begin{aligned}
\frac{N}{ L} \to \rho 
\  ,\quad  
dL \to \infty\  \hbox{and}\quad
\frac{dL}{N}\to 0\,.
\end{aligned}
\end{equation*}
Then
    for any $\mu \in \mathcal{M}_1{(\overline{\R}_+)}$ there exist $\eta^{(L,N)}\in \Omega_{L,N}$ such that 
    \begin{align}
        \Hat{\mu}_{\#} {\eta}^{(L,N)} \weakconv \mu\,.
    \end{align}
\end{lemma}
\begin{proof}
We will see that it suffices to approximate discrete measures of the form 
\begin{align}\label{eq_discrete_prob_meas}
\alpha_{0}\delta_0 + \alpha_{\infty} \delta_{\infty}+
    \sum_{i=1}^m \alpha_i \delta_{p_i}\in \mathcal{M}_1(\overline{\R}_+)\,, 
\end{align}
with $p_i\in (0,\infty)$, $1\leq i \leq m$. Let $\nu$ be such a probability
measure.

We explicitly construct configurations in $\Omega_{L,N}$ that converge to $\nu$
under the map
\begin{equation*}
\begin{aligned}
 \Hat{\mu}^{(\eta)}=\Hat{\mu}^{(\eta)}_{L,N}
    = \sum_{x=1}^L \frac{\eta_x}{N} \delta_{dL\frac{\eta_x}{N}} 
\,,
\end{aligned}
\end{equation*}
cf. \eqref{eq_def_map_mv_dLinf}, when
considering the thermodynamic limit $N/L\to \rho \in [0,\infty]$.

First, we consider the point masses lying in $(0,\infty)$. For convenience, let
us introduce 
\begin{equation*}
\begin{aligned}
k_i:= \left\lfloor \frac{ N}{dL} p_{i}  \right\rfloor 
\quad \hbox{and} \quad
\#_i:= \left\lfloor \frac{ \alpha_i N }{k_i}
\right\rfloor \,.
\end{aligned}
\end{equation*}
Note that $k_i \to \infty$, as $N/L \to \rho$, by assumption.

Now, let $ \eta'$ be the vector given by gluing together vectors $(k_i,\ldots ,
k_i) \in \N^{\#_i}$, $1\leq i \leq m$. Note that 
\begin{equation*}
\begin{aligned}
\frac{1}{L} 
\sum_{i=1}^{m} \#_{i} 
\leqslant 
\sum_{i =1}^{m}  \frac{ \alpha_i\, N }{\frac{1}{2}k_i \, L}
\leqslant 
2d
\sum_{i =1}^{m}  \frac{ \alpha_i }{ p_{i}}
\end{aligned}\,,
\end{equation*}
and the r.h.s. vanishes because the sum is finite, recall $p_i >0$.
Hence, without loss of generality, we assume
that the resulting vector $\eta'$ lies in $\N^{L}$, 
as we can append zeros to the constructed vector until $ \eta '$ has length $L$. 
In particular, the relative number of empty sites converges to one.

It only remains to distribute the remaining particles to create the point
masses at zero and infinity. Thus, far only 
$ \#_{\Sigma}=\sum_{i=1}^{m} \#_{i} $ sites are occupied in $ \eta '$.
 We start by adding $k_\infty := \lfloor
\alpha_{\infty} N \rfloor $ particles onto the $\#_{ \Sigma}+1$-th
position of $ \eta'$ (which is empty). This corresponds to the point mass at
infinity.

Now, the number of allocated particles to $ \eta' $ is upper bounded by 
\begin{equation*} 
\begin{aligned}
k_{\infty}+
\sum_{i=1}^{m} k_i \, \#_{i} 
\leqslant 
\alpha_{\infty}N +
\sum_{i=1}^{m} \alpha_{i} N = (1- \alpha_{0}) N\,.
\end{aligned}
\end{equation*}
We distribute the remaining 
$k_0 := N- k_{\infty}-\sum_{i=1}^{m} k_i \, \#_{i}$ 
particles as uniform as possible between all empty sites of $ \eta'$ (there
are $ \#_{0}:=L- \#_{\Sigma}-1$ many). This yields the configuration
\begin{equation*} 
\begin{aligned} 
\eta_x:=
\begin{cases}
\eta'_x & \quad \text{if $1\leq x \leq \#_{\Sigma}+1$}\\
 \left\lfloor\frac{k_0}{\#_0} \right\rfloor + \mathds{1}_{x \in
\{\#_{\Sigma}+2,\ldots, 
\#_{\Sigma}+1+ (k_0 \mod \#_{0})\}}&
 \quad \text{otherwise}.
\end{cases}
\end{aligned}
\end{equation*}
Because the number of non-empty sites in $ \eta'$ was relatively vanishing, the
particles distributed on previously empty sites will
correspond to the point mass at zero. Note that $ \eta \in \Omega_{L,N}$ by
construction.

Indeed, the constructed particle configuration $ \eta$ approximates the
discrete measure arbitrarily well in the thermodynamic limit, as 
for every $f\in C_b(\Bar{\R}_+)$ we have 
\begin{align}
\Hat{\mu}^{(\eta)}(f) 
=
\sum_{x = \#_{\Sigma}+2}^{L} \frac{\eta_{x}}{N} f \left(  \frac{dL}{N}
\eta_{x}   \right) 
+
\frac{k_{\infty}}{N} 
f \left(  \frac{dL}{N} k_{\infty}  \right)
+
\sum_{i=1}^{m} \frac{k_i}{N} \#_{i} \,
f \left( \frac{dL}{N} k_{i} \right)\,,
\end{align}
which converges to
\begin{equation*}
\begin{aligned}
\lim_{N/L \to \rho} \Hat{\mu}^{(\eta)}(f) 
=
\alpha_{0}f(0)
+
\alpha_{\infty}f(\infty)
+
\sum_{i=1}^{m} 
\alpha_{i} f(p_{i} ) = \nu(f)\,.
\end{aligned}
\end{equation*}
Now, for every $ \mu \in \mathcal{M}_{1}(\Bar{\R}_+)$ there exists a sequence $(\nu_n)_{n\in\N}$ of measures of the
form \eqref{eq_discrete_prob_meas} such that $\nu_n\weakconv \mu$. 
Moreover, each $\nu_{n}$ can be approximated by a sequence $(\eta_n^{(L,N)})_{L,N}$ following the above approach. 
Hence, we can construct a sequence of configurations with the desired property using a diagonal argument. 
For the sake of clarity we write out the details w.r.t the thermodynamic limit explicitly: we consider $N_j,L_j\to \infty$ such that $\lim_{j\to \infty}N_j/L_j= \rho>0$. Then for every $j$ we choose $n_j>n_{j-1}$ such that 
\begin{align} 
    d(\nu_{n_j}, \Hat{\mu}_{\#} \eta_{n_j}^{(L_j,N_j)} )\leq 2^{-j} 
\end{align}
where $d(\cdot, \cdot)$ denotes an appropriate metric, e.g. the L\'evy-Prokhorov metric. Thus, 
\begin{align} 
    d(\mu, \Hat{\mu}_{\#} \eta_{n_j}^{(L_j,N_j)} )
    \leq d(\mu,\nu_{n_j}) +2^{-j} \to 0\,, \quad \text{as }  j\to \infty\,,
\end{align}
which completes the proof.
\end{proof}


\section{Convergence to a Fleming-Viot process}\label{app_FV}

In this appendix we outline briefly how to prove convergence of the inclusion process to a
Fleming-Viot process with with mutation operator $A_{FV}$ \eqref{eq:FV_mut}, recall
\begin{equation*}
\begin{aligned}
        A_{FV}h(u) = \theta \int_0^1 [h(v)-h(u)]\,dv\,.
\end{aligned}
\end{equation*}
The generator of the Fleming-Viot process, when applied to a cylindrical test function of the form $H(\nu)=\nu(h_1) \cdots \nu(h_n)$, is given by
\begin{align}\label{eq_def_FV}
\mathcal{L}_{FV} H(\nu) &=2\sum_{1\leq k<l\leq n} \big(\nu (h_k h_l )-\nu (h_k
)\nu (h_l )\big) \prod_{j\neq k,l} \nu (h_j )\\
&\qquad +\sum_{1\leq k\leq n} \nu (A_{FV}
h_k )\prod_{j\neq k} \nu (h_j )\,. \nonumber
\end{align}
Instead of considering the generator $\mathfrak{L}_{L,N}$
\eqref{eq_def_ip_dynamics} on $\Omega_{L,N}$, it is more convenient to define the inclusion process on a
state space keeping track of particle positions. More precisely, 
we define
\begin{equation*}
\begin{aligned}
S_{L,N}:= \Lambda^{N}_{L}\,, \quad \text{with }\quad \Lambda_{L}:= \{1, \ldots,
L\}\,.
\end{aligned}
\end{equation*}
Now, the (labelled) inclusion process is described by the infinitesimal
generator 
\begin{align}\label{eq_gen_label_IP}
    \mathfrak{G}_{L,N}g(\sigma) 
    =
    \sum_{i,j=1}^N [g(\sigma^{i\to \sigma_j})- g(\sigma)]
    +
    d
    \sum_{i=1}^N
    \sum_{x\in \Lambda}
    [g(\sigma^{i\to x})- g(\sigma)]\,,
\end{align}
where
\begin{align}
    \sigma^{i\to x}_j
    :=
    \begin{cases}
    x& \quad \text{if } i=j\,,\\
    \sigma_j &\quad \text{if } i\neq j\,,
    \end{cases}
\end{align}
denotes the updated position after the $i$-th particle jumped onto site $x$.
The generator $\mathfrak{G}_{L,N}$ characterises a Markov process on $S_{L,N}$
which we denote by $(\sigma^{(L,N)} (t))_{ t \geqslant 0}$.
 We can recover the corresponding unlabelled particle configuration by
\begin{equation*}
\begin{aligned}
\iota :S_{L,N} \to\Omega_{L,N}\quad\mbox{with}\quad \iota ( \sigma)_x:=
\sum_{i =1}^{N} \mathds{1}_{\sigma_{i} = x} \,,
\end{aligned}
\end{equation*}
and in particular we have $
\mathfrak{L}_{L,N}f (\iota(\sigma))=\mathfrak{G}_{L,N}f(\iota (\cdot))( \sigma)
$.

Again, we interpret particle configurations as probability
measures on $[0,1]$. However, now we consider the embedding 
\begin{equation}\label{e:nu_embedding}
\begin{aligned}
\nu^{(\cdot)}:\, 
\sigma \mapsto \frac{1}{N} \sum_{i =1}^{N} \delta_{\frac{\sigma_{i}}{L}} \in
\mathcal{M}_{1} ([0,1]) \,,
\end{aligned}
\end{equation}
where rescaled particle
locations are encoded on the `type space' $[0,1]$.
Now the convergence of processes under the embedding \eqref{e:nu_embedding} follows from approximation of generators in analogy to our main result. 
Again, we start with test functions of the form 
\begin{equation*}
\begin{aligned}
\nu^{( \sigma)}(h) =\frac{1}{N} \sum_{i =1}^{N} h \left( \tfrac{ \sigma_{i}}{L}  \right)
\,, \quad h \in C^{3}([0,1])\,,
\end{aligned}
\end{equation*}
in which case the action of $\mathfrak{G}_{L,N}$ reads
\begin{equation*}
\begin{aligned}
\mathfrak{G}_{L,N}\nu^{( \sigma)}(h) 
=
  \sum_{i,j=1}^N \frac{1}{N} [h(\tfrac{\sigma_j}{L})- h(\tfrac{\sigma_j}{L})]
    +
    d\sum_{i=1}^N\sum_{x\in \Lambda} \frac{1}{N}
    [h(\tfrac{x}{L})- h(\tfrac{\sigma_i}{L})]\,.
\end{aligned}
\end{equation*}
Because the first sum on the r.h.s. vanishes by symmetry, we
are only left with 
\begin{equation*}
\begin{aligned}
\mathfrak{G}_{L,N} \nu^{( \sigma)}(h) 
=
d\sum_{i=1}^N\sum_{x\in \Lambda} \frac{1}{N}
    [h(\tfrac{x}{L})- h(\tfrac{\sigma_i}{L})]
 =
    dL\, \cdot \,  \nu^{( \sigma)} \left( \frac{1}{L} \sum_{x\in \Lambda}
    h(\tfrac{x}{L})- h \right)\,,
\end{aligned}
\end{equation*}
which implies the uniform convergence 
\begin{equation*}
\begin{aligned}
\lim_{N/L \to \rho} 
\sup_{\sigma \in \Lambda^{N}} 
\left|\mathfrak{G}_{L,N}(\nu^{( \cdot )}(h))(\sigma) - \nu^{( \sigma)}(A_{FV}h)  \right|=0\,.
\end{aligned}
\end{equation*}
By considering cylindrical test-functions, this convergence can be extended to
a core of the Fleming-Viot process with generator $\mathcal{L}_{FV}$
\eqref{eq_def_FV} in full analogy to our main results in
Section~\ref{sec_dLfin}. We leave out further details. 

Overall, this yields convergence of the (labelled) inclusion process in the
following sense: if $ \nu_{\#} \sigma^{(L,N)}(0) \weakconv \nu_{0}$ then 
\begin{equation*}
\begin{aligned}
\left(\nu_{\#} \sigma^{(L,N)} (t)\right)_{t \geqslant 0}
\weakconv
( \nu_{t})_{t \geqslant 0}\,, \quad \text{ in } D([0, \infty),
\mathcal{M}_{1}([0,1]))\,,
\end{aligned}
\end{equation*}
where $( \nu_{t})_{t \geqslant 0}$ denotes the Fleming-Viot process generated
by \eqref{eq_def_FV} with initial condition $ \nu_{0}$.